\documentclass[11pt]{article}

\usepackage{amsmath}
\usepackage{amssymb}
\usepackage{theorem}
\usepackage{euscript}
\usepackage{amscd}
\usepackage[arrow,curve,matrix,tips,frame]{xy}
\usepackage{pstricks}
\usepackage{pstcol,pst-node,pst-tree}

\renewcommand{\title}[1]{
     \addvspace{3\baselineskip}  
     \begin{center} \LARGE \bf #1
     \end{center}
     \addvspace{2\baselineskip}}   

\renewcommand{\author}[1]{
     \addvspace{-1\baselineskip}  
     \begin{center} \large \sc #1
     \end{center}
     \addvspace{2\baselineskip}}

\makeatletter

\def\section{%
        \@startsection{section}{1}{\z@}%
        {8ex plus 6ex minus 3ex}{\baselineskip}%
        {\normalfont\large\scshape\centering}%
        }

\renewcommand{\paragraph}[1]{{\par\removelastskip\vskip.5\baselineskip%   
         \indent{\itshape{#1}}{\ifperiod.\else\global\periodtrue\fi}%
         \rm \ignorespaces}}

\let\goth=\mathfrak
\let\calligraphy=\mathcal

\def\CC{{\mathbb C}}

\def\PP{{\mathbb P}}
\def\QQ{{\mathbb Q}}

\def\ZZ{{\mathbb Z}}

\def\Ff{{\calligraphy F}}
\def\Gg{{\calligraphy G}}
\def\Hh{{\calligraphy H}}
\def\Ii{{\calligraphy I}}
\def\Jj{{\calligraphy J}}

\def\Ll{{\calligraphy L}}

\def\Oo{{\calligraphy O}}

\def\Zz{{\calligraphy Z}}

\def\AAA{{\goth A}}

\def\FFF{{\goth F}}

\def\VVV{{\goth V}}

\def\hbar{{\,\overline{\!h}}}

\def\wbar{{\,\overline{\!w}}}

\def\Ktilde{{\,\widetilde{\!K}}}

\def\Stilde{{\,\widetilde{\!S}}}

\def\fhat{{\,\hat{\!f}}}
\def\Ghat{{\,\widehat{\!G}}}

\def\Hhat{{\,\widehat{\!H}}}

\def\oomega{{\boldsymbol{\omega}}}

\def\cc{{\boldsymbol{c}}}

\def\eE{{\boldsymbol{E}}}

\def\fF{{\boldsymbol{F}}}

\def\ww{{\boldsymbol{w}}}

\def\Bl{\operatorname{Bl}}
\def\card{\operatorname{card}}

\def\Hom{\operatorname{Hom}}
\def\Homcal{\Hh om}
\def\Im{\operatorname{Im}}

\def\ord{\operatorname{ord}}

\def\Res{\operatorname{Res}}

\def\Specbf{\mathbf{Spec\,}}
\def\Sing{\operatorname{Sing}}

\def\Trace{\operatorname{Tr}}

\let\lra=\longrightarrow
\newcommand{\arrow}[1]{\stackrel{#1}{\longrightarrow}}

\def\ie{{\it i.e.}~}

\def\inv{^{-1}}

\let\phi=\varphi
\let\epsilon=\varepsilon

\newcommand{\floor}[1]{\left\lfloor#1\right\rfloor}
\newcommand{\ceil}[1]{\left\lceil#1\right\rceil}

\def\theoname{Theorem}
\def\lemmaname{Lemma}
\def\propositionname{Proposition}
\def\notationname{Notation}
\def\corollaryname{Corollary}
\def\conjecturename{Conjecture}
\def\remarkname{Remark}
\def\remarksname{Remarks}
\def\examplename{Example}
\def\examplesname{Examples}
\def\definitionname{Definition}
\def\definitionsname{Definitions}
\def\notationname{Notation}

\def\proofname{Proof}

\def\Dquad{\hskip 0.6em plus .02em minus .2em}  
\def\Dpar{\belowdisplayskip=0pt\belowdisplayshortskip=0pt\par}

\def\bigpenalty{\interlinepenalty=\@M}
\def\smallpenalty{\interlinepenalty=100}

\newif\ifperiod \periodtrue 

\def\D@makemargins{%
  \labelsep=0pt
  \itemindent=0pt
  \labelwidth=0pt}

\def\D@restoremargins{%
  \labelsep=5pt
  \itemindent=0pt
  \leftmargin=5mm  
  \labelwidth=\leftmargin \advance\labelwidth by -\labelsep}

\def\th@Dindent{\hspace\parindent}
\def\th@Dheadingshape{\scshape}

\gdef\th@DthAndSuchtheo{%
  \D@makemargins%
  \def\@begintheorem##1##2{%
  \item[]\th@Dindent{\th@Dheadingshape ##1~\rm ##2.}\Dquad         
        \D@restoremargins}%
  \def\@opargbegintheorem##1##2##3{\def\next{##3}%
  \item[]\th@Dindent{\th@Dheadingshape ##1~\rm ##2\ifx\next\empty
  \else\ {\normalfont(##3)}\fi.}         
        \D@restoremargins}}

\gdef\th@DthAndSuchtheostar{%
  \D@makemargins%
  \def\@begintheorem##1##2{%
  \item[]\th@Dindent{\th@Dheadingshape ##1.}\Dquad     
        \D@restoremargins}%
  \def\@opargbegintheorem##1##2##3{\def\next{##3}%
  \item[]\th@Dindent{\th@Dheadingshape ##1\ifx\next\empty
  \else\ ##3\fi.}\Dquad         
        \D@restoremargins}}

\gdef\th@DthAndSuchliketheo{
  \D@makemargins%
  \def\@begintheorem##1##2{%
    \@latex@error{likethm: You must provide an argument in square brackets,
    though it may be empty [] !}%
    }%
  \def\@opargbegintheorem##1##2##3{%
        \def\next{##3}\ifx\next\empty\item[\th@Dindent]\else
        \item[]\th@Dindent{\th@Dheadingshape \next.}\Dquad\fi
        \D@restoremargins}}

\gdef\th@DdefAndSuch{%
  \D@makemargins%
  \def\@begintheorem##1##2{%
  \item[]\th@Dindent{\def@Dheadingshape ##1~\rm ##2.}\Dquad         
        \D@restoremargins}%
  \def\@opargbegintheorem##1##2##3{\def\next{##3}%
  \item[]\th@Dindent{\def@Dheadingshape ##1~\rm ##2\ifx\next\empty
  \else\ {\normalfont(##3)}\fi.}         
        \D@restoremargins}}

\gdef\th@DdefAndSuchStar{%
  \D@makemargins%
  \def\@begintheorem##1##2{%
  \item[]\th@Dindent{\def@Dheadingshape ##1.}\Dquad     
        \D@restoremargins}%
  \def\@opargbegintheorem##1##2##3{\def\next{##3}%
  \item[]\th@Dindent{\th@Dheadingshape ##1\ifx\next\empty
  \else\ ##3\fi.}\Dquad         
        \D@restoremargins}}

\def\th@Dheadingshape{\scshape}
\def\def@Dheadingshape{\itshape}

\theoremstyle{DthAndSuchliketheo}
\theorembodyfont{\bigpenalty\itshape}   
\newtheorem{likethm}{}
\theorembodyfont{\rmfamily}  

\theoremstyle{DthAndSuchtheostar}
\theorembodyfont{\bigpenalty\itshape}
\newtheorem{thm*}{\theoname}
\newtheorem{lem*}{\lemmaname}
\newtheorem{pro*}{\propositionname}
\newtheorem{cor*}{\corollaryname}
\newtheorem{conjecture*}{\conjecturename}
\theorembodyfont{\smallpenalty\rmfamily}
\newtheorem{notation*}{\notationname}
\newtheorem{exa*}{\examplename}
\newtheorem{examples*}{\examplesname}
\theoremstyle{DdefAndSuchStar}
\theorembodyfont{\smallpenalty\rmfamily}
\newtheorem{definition*}{\definitionname}
\newtheorem{definitions*}{\definitionsname}
\newtheorem{rem*}{\remarkname}
\newtheorem{remarks*}{\remarksname}
\theoremstyle{DthAndSuchtheo}
\theorembodyfont{\bigpenalty\itshape}
\newtheorem{thm}{\theoname}[section]
\newtheorem{lem}[thm]{\lemmaname}
\newtheorem{pro}[thm]{\propositionname}
\newtheorem{cor}[thm]{\corollaryname}

\theorembodyfont{\smallpenalty\rmfamily}
\newtheorem{exa}[thm]{\examplename}

\theoremstyle{DdefAndSuch}
\theorembodyfont{\smallpenalty\rmfamily}
\newtheorem{definition}[thm]{\definitionname}
\newtheorem{rem}[thm]{\remarkname}

\newtheorem{remarks}[thm]{\remarksname}

\theoremstyle{DthAndSuchtheo}               
\theorembodyfont{\itshape}

\newcommand{\proof}[1][]{{\par\removelastskip\vskip.6\baselineskip   
    \noindent\th@Dindent\def\next{#1}%
    {\itshape\proofname\ifx\next\empty\else\next\fi\ifperiod.%
      \else\global\periodtrue\fi\Dquad}%
    \clubpenalty=5000\rm\ignorespaces}\setcounter{step}{0}}
\newcounter{step}
\newcommand{\step}{\stepcounter{step}
  \par\indent{\itshape Step \thestep .\hspace{1ex}}}

\newcommand{\likeproof}[1][]{{\par\removelastskip\vskip.6\baselineskip
    \noindent\th@Dindent\def\next{#1}%
    {\itshape\ifx\next\empty\else\next\fi\ifperiod.%
      \else\global\periodtrue\fi\Dquad}%
    \clubpenalty=5000\rm\ignorespaces}\hspace{-2pt}\setcounter{step}{0}}

\def\qed{{\ifmmode\hskip 6mm plus 1mm minus 3mm{$\square$}
    \else
    \hfil\penalty50\hskip1em\null\nobreak\hfil
    {\hfill $\square$\parfillskip=0pt\finalhyphendemerits=0
      \let\par=\endgraf\par}
    \fi
    \Dpar\penalty-150\vskip.6\normalbaselineskip}}

\makeatother

\newpsobject{showgrid}{psgrid}{subgriddiv=1,griddots=10,gridlabels=6pt}

\begin{document}

\title{The irregularity of cyclic multiple planes after Zariski}
\author{Daniel Naie}

\begin{flushleft}
\end{flushleft}

\begin{abstract}
A formula for the irregularity of a cyclic multiple plane associated
to a branch curve that has arbitrary singularities and is transverse
to the line at infinity is established.  The irregularity is expressed
as a sum of superabundances of linear systems associated to some
multiplier ideals of the branch curve and the proof rests on the
theory of standard cyclic coverings.  Explicit computations of
multiplier ideals are performed and some applications are presented.
\end{abstract}

\section*{Introduction}

Let $f(x,y)=0$ be an affine equation of a curve $B\subset\PP^2$ and
$H_\infty$ be the line at infinity.  The projective surface
$S_0\subset\PP^3$ defined by the affine equation $z^n=f(x,y)$ is
called the $n$-cyclic multiple plane associated to $B$ and $H_\infty$.
In \cite{Za}, Zariski obtains his famous result for the irregularity
of certain $n$-cyclic multiple planes.

\begin{likethm}[Zariski's Theorem]
Let $B$ be an irreducible curve of degree $b$, transverse to the line
at infinity $H_\infty$ and with only nodes and cusps as singularities.
Let $S_0\subset\PP^3$ be the $n$-cyclic multiple plane associated to
$B$ and $H_\infty$, and let $S$ be a desingularization of $S_0$.  The
surface $S$ is irregular if and only if $n$ and $b$ are both divisible
by $6$ and the linear system of curves of degree $5b/6-3$ passing
through the cusps of $B$ is superabundant.  In this case,
\[
  q(S) = h^1(\PP^2,\Ii_\Zz(-3+\frac{5b}{6})),
\]
where $\Zz$ is the support of the set of cusps.
\end{likethm}

The aim of this paper is to present a generalization of Zariski's
Theorem to a branch curve that has arbitrary singularities and is
transverse to the line at infinity bringing to the fore the theory of
cyclic coverings as developed in \cite{Pa}.  The irregularity will be
expressed as a sum of superabundances of linear systems associated to
some multiplier ideals of the branch curve $B$.  We refer to
\cite{EiLa} for the notion of multiplier ideal.  To state the main
result in Section \ref{s:theIrregularity}, we recall here that if the
rational $c$ varies from a very small positive value to $1$, then one
can attach a collection of multiplier ideals $\Jj(cB)$ that starts at
$\Oo_{\PP^2}$, diminishes exactly when $c$ equals a jumping
number---they represent an increasing discrete sequence of
rationals---and finally ends at $\Ii_B$.

\begin{thm*}[(\ref{th:q})]
Let $B$ be a plane curve of degree $b$ and let $H_\infty$ be a line
transverse to $B$.  Let $S$ be a desingularization of the $n$-cyclic
multiple plane associated to $B$ and $H_\infty$.  If
$J(B,n)$ is the subset of subunitary jumping numbers of $B$ that live
in $\frac{1}{\gcd(b,n)}\ZZ$, then
\[
  q(S) = \sum_{\xi\in J(B,n)}h^1(\PP^2,\Ii_{Z(\xi B)}(-3+\xi b)),
\]
where $Z(\xi B)$ is the subscheme defined by the multiplier ideal
$\Jj(\xi B)$. 
\end{thm*}

In case the singularities of $B$ are locally given by $x^p=y^q$ such
as equations, explicit computations of the jumping numbers and of the
multiplier ideals will enable us to apply the above theorem to various
examples in Section \ref{s:applications}.  An example in Remark
\ref{r:nonTransverse} shows that the irregularity may jump in case the
position of $H_\infty$ with respect to $B$ becomes special. 

Generalizations of Zariski's Theorem are discussed in several papers
and the proofs are based on different points of view.  First,
Zariski's original argument divides naturally into three parts.  He
describes the canonical system of $S$ in terms of the conditions
imposed by the singularities of $S_0$ that correspond to the cusps.
Then he establishes the formula
\begin{equation}
  \label{eq:qZInGeneral}
  q(S)
  = \sum_{k=n-\floor{n/6}}^{n-1}
    h^1(\PP^2,\Ii_\Zz(-3+\ceil{\frac{kb}{n}})),
\end{equation}
where $\Zz$ denotes the support of the set of cusps.  To finish, he
invokes the topological result proved in \cite{Za2}: {\it If $n$ is
the power of a prime and $B$ is irreducible, then the $n$-cyclic
multiple plane is regular.}  The theorem follows from the examination
of the terms that vanish in the previous sum when the degree of the
cyclic multiple plane covers an unbounded sequence of powers of
primes.

Second, in \cite{Es}, Esnault establishes a formula, similar to
(\ref{eq:qZInGeneral}), for the irregularity of the $b$-cyclic
multiple plane $S_0$, where $b$ is the degree of the branch curve $B$
that possesses arbitrary isolated singularities.  She uses the
techniques of logarithmic differential complexes, the existence of a
mixed Hodge structure on the complex cohomology of the associated
Milnor fibre---the complement of $S_0$ with respect to the plane that
contains $B$---and Kawamata-Viehweg Vanishing Theorem.  In \cite{Ar},
Artal-Bartolo interprets Esnault's formula for irregularity and
applies it to produce two new Zariski pairs.  Two plane curves
$B_1,B_2\subset\PP^2$ are called a Zariski pair if they have the same
degree and homeomorphic tubular neighbourhood in $\PP^2$, but the
pairs $(\PP^2,B_1)$ and $(\PP^2,B_2)$ are not homeomorphic.  Zariski
was the first to discover that there are two types of plane sextics
with six cuspidal singularities: there are the ones where the cusps
lie on a conic, and the ones where the cusps don't lie on a plane
conic.  In \cite{Va}, Vaqui\'e gives a formula for the irregularity of
a cyclic covering of degree $n$ of a nonsingular algebraic surface $Y$
ramified along a reduced curve $B$ of degree $b$ with respect to some
projective embedding and a nonsingular hyperplane section $H$ that
intersects $B$ transversally.  His formula is stated in terms of
superabundances of the set of singularities of $B$ and the proof also
uses the techniques of logarithmic differential complexes.  The
superabundances involved are given by ideal sheaves that coincide in
fact to the multiplier ideals.  Vaqui\'e's paper is one among several
to introduce the notion of multiplier ideals implicitly and we refer
to \cite{EiLa} for this issue.

Third, in \cite{Li1}, Libgober applies methods from knot theory to
study the $n$-multiple plane $S_0$.  His results are expressed in
terms of Alexander polynomials and extend Zariski's Theorem to
irreducible curves $B$ with arbitrary singularities and to lines
$H_\infty$ with arbitrary position with respect to $B$.  
Later on, in \cite{Li2,Li3,Li4}, he deals with the case of reducible
curves $B$ having transverse intersection with the line at infinity
and the irregularity of the multiple plane is expressed using
quasiadjunction ideals.  The technique is based on mixed Hodge theory,
and the result is a particular case in a vaster study, pursued in the
above mentioned papers, where the homotopy groups of the complements
of various divisors in smooth projective varieties are explored.
These groups are related to the Hodge numbers of cyclic or more
generally abelian coverings ramified along the considered divisors, as
well as to the position of their singularities.  We refer the reader
to \cite{Li6} for more ample details and references and to \cite{Li5}
for the relation between the quasiadjunction ideals and the multiplier
ideals.

Our argument will follow Zariski's ideas.  The multiple plane is
transformed into a standard cyclic covering of the plane through a
sequence of blowing-ups of $\PP^3$.  Then an analog of the formula
(\ref{eq:qZInGeneral}) is obtained thanks to the theory of cyclic
coverings:
\begin{equation*}
  \label{eq:qInGeneral}
  q(S) 
  = \sum_{k=1}^{n-1}
  h^1(\PP^2,\Ii_{Z(\frac{k}{n}B)}(-3+\ceil{\frac{kb}{n}})).
\end{equation*}
Finally Theorem \ref{th:q} is established using Kawamata-Viehweg-Nadel
Vanishing Theorem.

\begin{rem*}
The above formula coincides with Vaqui\'e's in \cite{Va} when the
latter is interpreted for a plane curve $B$ and a line $H$ transverse
to it.  At the same time, Vaqui\'e's formula in its general form might
be obtained by the argument we make use of in establishing Theorem
\ref{th:q} if Vaqui\'e's general setting were to be considered.
\end{rem*}

The paper is organized as follows.  In Section \ref{s:preliminaries}
the theory of cyclic coverings and some facts about multiplier ideals
are recalled.  Next, in Section \ref{s:theIrregularity} it is shown
how through a sequence of blowing-ups a cyclic multiple plane is
transformed into a standard cyclic covering of the plane and Theorem
\ref{th:q} is proved.  Explicit computations of the jumping numbers
and multiplier ideals are performed in Section \ref{s:clusters}, using
the theory of clusters.  Finally, in Section \ref{s:applications} some
applications are presented.

\paragraph{Notation and conventions}
All varieties are assumed to be defined over $\CC$.  Standard symbols
and notation in algebraic geometry will be freely used.  Moreover, if
$D$ is a divisor on the variety $X$, we shall often write $H^i(X,D)$
and $h^i(X,D)$ instead of $H^i(X,\Oo_X(D))$ and $h^i(X,\Oo_X(D))$
respectively.  If $\Ll$ is an invertible sheaf on $X$, then we shall
regularly denote by $L$ a divisor such that $\Ll\simeq\Oo_X(L)$.

\paragraph{Acknowledgements} I started this paper during a one week
stay at the University of Pisa in the spring of 2004.  I would like to
thank Rita Pardini for her hospitality and for the friendly talks we
had.

The paper owes Mihnea Popa its present form.  I would like to record
my debt to his reading of a preliminary version in the autumn of 2005
and to his encouragements to generalize the results I obtained at that
time.

Finally, I would like to thank my colleagues Laurent Evain and
Jean-Philippe Monnier for the conversations they put up with
throughout this period.

\section{Preliminaries}
  \label{s:preliminaries}
We shall summarize, in a form convenient for further use, some
properties of cyclic coverings and of multiplier ideals.

\subsection{Cyclic coverings} 
Let $X$ be a variety and let $G$ be the finite abelian group of order
$n$.  If $G$ acts faithfully on $X$, then the quotient $Y=X/G$ exists
and $X$ is called an abelian covering of $Y$ with group $G$.  The map
$\pi:X\to Y$ is a finite morphism, $\pi_\ast\Oo_X$ is a coherent sheaf
of $\Oo_Y$-algebras, and $X\simeq\Specbf (\pi_\ast\Oo_X)$.  

If $X$ is normal and $Y$ is smooth, then $\pi$ is flat which is
equivalent to $\pi_\ast\Oo_X$ locally free.  The action of $G$ on
$\pi_\ast\Oo_X$ decomposes it into the direct sum of eigensheaves
associated to the characters $\chi\in\Ghat$, 
\[
 \pi_\ast\Oo_X=\bigoplus_{\chi\in\Ghat}\Ll_\chi\inv.
\]
The action of $G$ on $\Ll_\chi$ is the multiplication by
$\chi$ and $\Ll_1=\Oo_Y$. 

To understand the ring structure of $\pi_\ast\Oo_X$ we suppose that
every component $D$ of the ramification locus is $1$-codimensional.
Such a component is associated to its stabilizer subgroup 
$H\subset G$ and to a character $\psi\in\Hhat$ that
generates $\Hhat$: $\psi$ corresponds to the induced representation of
$H$ on the cotangent space to $X$ at $D$.
Dualizing the inclusion $H\subset G$, such a couple $(H,\psi)$ is
equivalent to a group epimorphism $f:\Ghat\to\ZZ/m_f$, where
$m_f=|H|$; for any $\chi\in\Ghat$, the induced representation
$\chi|_H$ is given by $\psi^{f(\chi)^\bullet}$.  

Here and later on, $a^\bullet$ denotes the smallest non-negative
integer in the equivalence class of $a\in\ZZ/m$, and $\FFF$ the set of
all group epimorphisms from $\Ghat$ to different $\ZZ/m\ZZ$.

Let $B_f$ be the divisor whose components belong to the branch locus
and are exactly those covered by components of the ramification locus
associated to the group epimorphism $f$.  The ring structure is given
by the following isomorphisms (see \cite{Pa}): for any
$\chi,\chi'\in\Ghat$,  
\begin{equation} \label{eq:product}
  \Ll_\chi\otimes\Ll_{\chi'} \simeq 
  \Ll_{\chi\chi'}\otimes\bigotimes_{f\in\FFF}
  \Oo_Y(\epsilon(f,\chi,\chi')B_f)
\end{equation}
with $\epsilon(f,\chi,\chi')=0$ or $1$, depending on whether or not
$f(\chi)^\bullet+f(\chi')^\bullet<|\Im \fhat|$. 

The next proposition is formulated for cyclic groups, since it is in
this case that will be used in the sequel.  We refer again to
\cite{Pa} for the case of abelian groups.

\begin{pro}
Let $\pi:X\to Y$ be a cyclic covering with $X$ normal, $Y$ smooth and
every component of the ramification locus $1$-codimensional. 
If $\chi$ generates $\Ghat$, then for every $k=1,\ldots,n$,
\begin{equation}\label{eq:k}
  L_{\chi^k}
  \equiv kL_\chi -
  \sum_{f\in\FFF}\floor{\frac{kf(\chi)^\bullet}{m_f}}B_f.
\end{equation}
In particular, for $k=n$ equation (\ref{eq:k}) becomes
\begin{equation}\label{eq:n}
  nL_\chi
  \equiv\sum_{f\in\FFF}\,
  [G:\Im\fhat]\, f(\chi)^\bullet B_f. 
\end{equation}
\end{pro}
\proof
For the proof we need to define the sequence
$(\zeta_k^{m,r})_{k\geq0}$: for $m$ and $r$ fixed positive integers
with $r\leq m$, and for $k\geq0$, put $\zeta_k^{m,r}=1$ if
$\floor{kr}_m^\bullet<r$, and $\zeta_k^{m,r}=0$ otherwise.  Obviously
this sequence is $m$-periodic, $\zeta_m^{m,r}=\zeta_0^{m,r}=1$ and in
case $m>r$, $\zeta_1^{m,r}=0$.

Now, from the hypotheses, $\chi$ spans the group of characters.
Taking $\chi'=\chi^{j-1}$ in (\ref{eq:product}) we get 
\[
  L_\chi + L_{\chi^{j-1}} \equiv 
  L_{\chi^j} + \sum_{f\in\FFF}\zeta_j^{m_f,f(\chi)^\bullet}B_f,   
\]
since 
$\epsilon(f,\chi,\chi^{j-1})=1$,
$f(\chi)^\bullet+f(\chi^{j-1})^\bullet\geq m_f$, 
$(f(\chi)+f(\chi^{j-1}))^\bullet<r$ and 
$[jf(\chi)^\bullet]_{m_f}^\bullet<r$ are equivalent.
Then, summing over $j$ from $1$ to $k$,
\[
  L_{\chi^k} \equiv kL_\chi
  - \sum_{j=1}^k\sum_{f\in\FFF}\zeta_j^{m_f,f(\chi)^\bullet}B_f.
\]
But $\sum_{j=1}^k\zeta_j^{n,b}$ represents the number of $1$'s among
the first $k$ terms in the sequence $(\zeta_j^{n,b})_j$, hence
$\sum_{j=1}^k\zeta_j^{n,b}=\floor{kb/n}$ and (\ref{eq:k}) follows. 
Formula (\ref{eq:n}) is obvious, since $\chi^n=1$.
\qed

Conversely, to every set of data $\Ll_\chi$, $B_f$, with $f\in\FFF$,
that satisfies (\ref{eq:n}), using (\ref{eq:k}), we define the line
bundles $\Ll_{\chi^k}$  and associate in a natural way 
{\it the standard cyclic covering}
$\pi:\Specbf(\oplus_k\Ll_{\chi^k}\inv)\to Y$, unique up to
isomorphisms of cyclic coverings.
The line bundles $\Ll_{\chi^k}$ verify equation (\ref{eq:product})
and, consequently, $\oplus_k\Ll_{\chi^k}\inv$ is endowed with a ring
structure.

Now, the standard covering thus obtained may not be normal; in fact it
is not normal precisely above the multiple components of the branch
locus (see \cite{Pa} Corollary 3.1).

\subsection{The normalization procedure for standard cyclic coverings}
Let $f:\Ghat\to\ZZ/m_f$ be a group epimorphism, so
$m_f=\ord(\Im\fhat)$, and let $B_f=rC+R$, with $C$ irreducible and not
a component of $R$, and $r\geq 2$.  $X$ is not normal along the
pull-back of $C$.  The normalization procedure along this multiple 
component of the branch locus splits into three steps showing how
to end up with a  new covering $X'\to X\to Y$, with $X'$
normal along the pull-back of $C$ (see \cite{Pa}).  The steps are
given by the comparison between the multiplicity $r$ and the order
$m_f$ of the stabilizer subgroup. 

\step
If $B_f=rC+R$ with $r\geq m_f$, then set $q$ and $r'$ by the Euclidean
division $r=qm_f+r'$, and construct a new set of building data by
putting 
\[
  L'_\chi \equiv L_\chi-f(\chi)^\bullet q\,C ,
  \quad  B'_f \equiv r'C+R 
  \quad\text{and}\quad
  B'_g \equiv B_g \quad\text{if }g\neq f.
\]

\step
If $B_f=rC+R$ with $r<m_f$ and $(r,m_f)=d>1$, then the natural
composition is considered
\[
  f':\Ghat\arrow{f}\ZZ/m_f\lra\ZZ/\frac{m_f}{d}.
\]
The integers $f(\chi)^\bullet$ and $f'(\chi)^\bullet$ are linked by
the relation  
$f(\chi)^\bullet=qm_f/d+f'(\chi)^\bullet$.
Put
\[
  L'_\chi \equiv L_\chi-q\frac{r}{d}\,C ,
  \quad B'_f \equiv R , 
  \quad B'_{f'} \equiv B_{f'}+\frac{r}{d}\,C
  \quad\text{and}\quad
  B'_g \equiv B_g \quad\text{if }g\neq f,f'
\]
in order to construct a `less non-normal' covering.  Notice that the
induced multiplicity and the corresponding subgroup order become
relatively prime. 

\step
If $B_f=rC+R$ with $r<m_f$ and $(r,m_f)=1$, then the composition 
\[
  f':\Ghat\arrow{f}\ZZ/m_f\arrow{r\cdot}\ZZ/m_f
\]
is considered.  As before, the integers $f(\chi)^\bullet$ and 
$f'(\chi)^\bullet$ are linked by
$r\cdot f(\chi)^\bullet=q m_f+f'(\chi)^\bullet$.
Put
\[
  L'_\chi \equiv L_\chi-q\,C ,
  \quad B'_f \equiv R ,
  \quad B'_{f'} \equiv B_{f'}+C 
  \quad\text{and}\quad
  B'_g \equiv B_g \quad\text{if }g\neq f,f'
\]
to get a new covering $X'$ and finish the normalization procedure
along $C$.

\begin{example}
On $\PP^2$ let $\Ll_\chi=\Oo(1)$ and 
$nL_\chi\equiv H_0+(n-1)H_\infty$, where $H_0$ and $H_\infty$ are two
fixed different lines.  Here the only functions $f:\Ghat\to\ZZ/n$
involved in (\ref{eq:product}) are given by $\chi\mapsto 1$ and
by $\chi\mapsto n-1$.  In this way, the standard $n$-cyclic covering
$S_0\to\PP^2$ is normal and has a singular point above $P$, the
intersection of $H_0$ and $H_\infty$.  To desingularize it, we
consider the blow-up surface $\Bl_P\PP^2$, with $E$ the exceptional
divisor and the induced cyclic covering $S\to\Bl_P\PP^2$.  We have
$nL_\chi\equiv H_0+(n-1)H_\infty+nE$ and the induced covering $S$ is
not normal above $E$. 
The normalization procedure leads to $S'\to\Bl_P\PP^2$ defined by
$nL'_\chi\equiv H_0+(n-1)H_\infty$, with $L'_\chi\equiv H-E$.  $S'$ is
a geometrically ruled surface and the pull-back of $E$ is a rational
section with self-intersection $-n$, \ie $S'$ is the Hirzebruch
surface $\fF_n$. 
\end{example}

\subsection{Multiplier ideals}\label{ss:multiplierIdeals}
Let $X$ be a smooth variety, $D\subset X$ be an effective
$\QQ$-divisor and $\mu:Y\to X$ be an embedded resolution for $D$.
Assume that the support of the $\QQ$-divisor $K_{Y|X}-\mu^\ast D$
is a union of irreducible smooth divisors with normal crossing
intersections.  Then $\mu_\ast\Oo_Y(K_{Y|X}-\floor{\mu^\ast D})$
is an ideal sheaf $\Jj(D)$ on $X$.  We will denote by $Z(D)$ the
subscheme defined by this ideal.  Hence  $\Ii_{Z(D)} = \Jj(D)$. 
Showing that $\Jj(D)$ is independent of the choice of the
resolution, see \cite{La}, we have:

\begin{definition*}
The ideal $\Jj(D)=\mu_\ast\Oo_Y(K_{Y|X}-\floor{\mu^\ast D})$ is
called the multiplier ideal of $D$.
\end{definition*}

The sheaf computing the multiplier ideal verifies the following
vanishing result: for every $i>0$,
$R^i\mu_\ast\Oo_Y(K_{Y|X}-\floor{\mu^\ast D})=0$.  
Therefore, applying the Leray spectral sequence, we obtain that for
every $i$  
\begin{equation}
  \label{eq:upDown}
  H^i(X, \Oo_X(K_X+L)\otimes\Ii_{Z(D)})
  = H^i(Y, 
  \Oo_Y(\mu^\ast K_X+\mu^\ast L + K_{Y|X}-\floor{\mu^\ast D})).
\end{equation}
Moreover,

\begin{likethm}[Kawamata-Viehweg-Nadel Vanishing Theorem]
Let $X$ be a smooth projective variety.  If $L$ is a Cartier divisor
and $D$ is an effective $\QQ$-divisor on $X$ such that $L-D$ is a nef
and big $\QQ$-divisor, then 
\[
  h^i(X,\Oo_X(K_X+L)\otimes\Ii_{Z(D)}) = 0
\] 
for every $i>0$.
\end{likethm}

\begin{likethm}[Definition-Lemma]\rm (see \cite{EiLa})
Let $B\subset X$ be an effective divisor and $P\in B$ be a fixed
point.  Then there is an increasing discrete sequence of rational
numbers $\xi_i:=\xi(B,P)$,
\[
  0 = \xi_0 < \xi_1 < \cdots
\]
such that 
\[
  \Jj(\xi B)_P = \Jj(\xi_i B)_P 
  \quad\text{for every}\quad
  \xi\in[\xi_i,\xi_{i+1}),
\]
and $\Jj(\xi_{i+1}B)_P\subset\Jj(\xi_i B)_P$.
The rational numbers $\xi_i$'s are called {\it the jumping
  numbers} of $B$ at $P$. 
\end{likethm}

\section{The irregularity of cyclic multiple planes}
  \label{s:theIrregularity}

\begin{thm}\label{th:q}
Let $B$ be a plane curve of degree $b$ and let $H_\infty$ be a line
transverse to $B$.  Let $S$ be a desingularization of the projective
$n$-cyclic multiple plane associated to $B$ and $H_\infty$.  If
$J(B,n)$ is the subset of subunitary jumping numbers of $B$ that live
in $\frac{1}{\gcd(b,n)}\ZZ$, then
\[
  q(S)
  = \sum_{\xi\in J(B,n)}h^1(\PP^2,\Ii_{Z(\xi B)}(-3+\xi b)),
\]
with $Z(\xi B)$ the subscheme defined by the multiplier ideal of 
$\xi B$.
\end{thm}

The proof splits naturally into four parts.  First, we show that there
is a sequence of blowing-ups of $\PP^3$ such that $S_1$, the
strict transform of the multiple plane $S_0\subset\PP^3$, becomes a
standard cyclic covering of the plane defined by 
$nL'_\chi\equiv B+(\beta n-b)H_\infty$, with 
$\beta=\ceil{b/n}$ and $\Ll'_\chi=\Oo_{\PP^2}(\beta)$.  Second we
choose a desingularization of $B$ such that its total transform on 
$\mu:Y\to\PP^2$ is a divisor with normal crossing intersections, a log
resolution.  It induces a standard cyclic covering $S_2$.  
We apply the normalization procedure to it and obtain a normalization
$S$ of $S_0$, defined by the line bundle $\Ll_\chi$ and that has only
Hirzebruch-Jung singularities. 
\[
\UseTips
 \newdir{ >}{!/-5pt/\dir{>}}
\xymatrix{   
  S \ar[rdd]_\pi \ar[rd] \\
  & S_2 \ar[r] \ar[d] & S_1 \ar[r] \ar[d] & S_0 \\
  & Y \ar[r]^\mu & \PP^2}
\]
Third, we compute the line bundles $\Ll_{\chi^k}$'s in terms of the
pull-back $\mu^\ast\Oo_{\PP^2}(1)$ and the exceptional configuration
on $Y$ and get the irregularity of $S$ as a sum of some $h^1$'s.
Finally, the result is obtained by applying the Kawamata-Viehweg-Nadel
Vanishing Theorem. 

The first step is given by:

\begin{pro}\label{p:multipleToCovering}
Let $S_0$ be the $n$-multiple plane associated to the curve $B$ of
degree $b$ and the line $H_\infty$.  There exists a sequence of
blowing-ups $S_1\to S_0$ such that $S_1$ is the standard
cyclic covering of the plane determined by 
\[
  nL'_\chi \equiv B+(\beta n-b)H_\infty,
\]
with $\beta=\ceil{b/n}$ and $L'_\chi\equiv\beta H$.
\end{pro}
\proof Let $[x_0,\ldots,x_3]$ be a homogeneous system of coordinates
in $\PP^3$.  Let $\Lambda$ be the plane defined by $x_3=0$, and let
$B,H_\infty\subset\Lambda$ be defined by $F(x_0,x_1,x_2)=0$ and
$x_0=0$, respectively.  $H_\infty$ will be called the line at
infinity.  The projective $n$-cyclic multiple plane $S_0\subset\PP^3$
is defined by $x_3^n = x_0^{n-b}F(x_0,x_1,x_2)$.

If $\deg B=b\leq n$, then things are easy.  The point $N$ of
homogeneous coordinates $[0,0,0,1]$ is not on $S_0$.  The complement
of the exceptional divisor $E$ in the blow-up of $\PP^3$ at $N$
coincides with the total space of the line bundle $\Oo_{\PP^2} (1)$.
Over any open subset $x_i\neq0$ of the projective plane $\Lambda$,
$i=0,1,2$, if $z=x_3/x_i$, then $z$ coincides with the tautological
section of $p^\ast \Oo_\Lambda(1)$ with $p:\Bl_N\PP^3-E\to\Lambda$.
The zero divisor of $p^\ast F-z$ defines $S_0$.  Hence $S_0$ is the
standard cyclic covering determined by
\[
  nL'_\chi \equiv B+(n-b)H_\infty,
\]
with $L'_\chi\equiv H$. 

If $\deg B=b>n$, then the situation is slightly more complicated since
now $N$ lies on $S_0$.  Let $\Xi$ be the plane spanned by $H_\infty$
and $N$.  In the open set $x_3\neq0$, $S_0$ is defined by
\[
  u_0^{b-n} = F(u_0,u_1,u_2),
\]
with $u_i=x_i/x_3$, $0\leq i\leq 2$.  First, we blow up the projective
space at $N$, $X_1=\Bl_N\PP^3\to\PP^3$.  $E\subset\Bl_N\PP^3$ denotes
again the exceptional divisor and $L_\infty\subset E$ the line that
correspond to $\Xi$, \ie $L_\infty=\Xi\cap E$.  The strict transform
$S_0$ is defined by
\[
  u_0^{(1)\,b-n} = u_1^{(1)\,n}F(u_0^{(1)},1,u_2^{(1)})
\]
on the subset $u_0=u_0^{(1)}u_1^{(1)}$, $u_1=u_1^{(1)}$, 
$u_2=u_1^{(1)}u_2^{(1)}$.  Notice that the line 
$L_\infty:u_0^{(1)}=u_1^{(1)}=0$ is contained in $S_0$. 
What we have to understand is the geometry of $S_0$ along 
$L_\infty$. 

Second, we see $X_1$ as 
$\PP(\Oo_\Lambda\otimes\Oo_\Lambda(1))\to\Lambda$ and make 
an elementary transform of $X_1$ along $L_\infty$.  We blow up 
$X_1$ along $L_\infty$ (the trace of the new exceptional divisor 
on $E$ is denoted by $L_\infty$).  Then, we contract the 
strict transform of\/ $\Xi$ to $L_\infty$.  We obtain 
$X_2=\PP(\Oo_\Lambda\otimes\Oo_\Lambda(2))\to\Lambda$.
The new exceptional divisor becomes an $\fF_1$ through $H_\infty$ 
and $L_\infty$, and will be denoted by $\Xi$.  
On $u_0^{(1)}=u_0^{(2)}$, $u_1^{(1)}=u_0^{(2)}u_1^{(2)}$, 
$u_2^{(1)}=u_2^{(2)}$, an equation for $S_0$ is 
\[
  u_0^{(2)\,b-2n} = u_1^{(2)\,n}F(u_0^{(2)},1,u_2^{(2)}),
\]
with $L_\infty:u_0^{(2)}=u_1^{(2)}=0$.

After $\beta-1$ elementary transforms along $L_\infty$, we get 
$S_0\subset\PP(\Oo_\Lambda\otimes\Oo_\Lambda(\beta))\to\Lambda$
with $S_0$ locally characterized by 
\[
  u_0^{(\beta)\,b-\beta n} 
  = u_1^{(\beta)\,n}F(u_0^{(\beta)},1,u_2^{(\beta)}).
\]
$E$ is defined by $u_0^{(\beta)}=0$ and $L_\infty$ by
$u_0^{(\beta)}=u_1^{(\beta)}=0$.  The new $\Xi$ is the Hirzebruch
surface $\fF_\beta$.  To finish, we put $z=1/u_1^{(\beta)}$ and look
at $x=u_0^{(\beta)}$ and $y=u_2^{(\beta)}$ as to local coordinates on
$\Lambda$.  Then 
\[
  S_0 : z^n=x^{\beta n-b}F(x,1,y).
\]
The complement of $E$ in $X_\beta$ seen through $p:X_\beta-E\to\Lambda$, 
coincides with the total space of $\Oo_\Lambda(\beta)$.  The
coordinate $z$ coincides with the tautological section of
$p^\ast\Oo_\Lambda(\beta)$.  We conclude that $S_0$ is the standard
cyclic covering determined by 
\[
  nL'_\chi \equiv B+(\beta n-b)H_\infty,
\]
with $L'_\chi\equiv\beta H$.
\qed

For the next step in the proof of Theorem \ref{th:q} we need several
preliminary results.

\begin{pro}\label{p:keyProposition}
Let $Y$ be smooth and let $\pi:X\to Y$ be a standard cyclic covering of
degree $n$ determined by  
\[
  nL_\chi\equiv\sum_{f\in\FFF}\,[G:\Im\fhat]\, f(\chi)^\bullet B_f.
\]
For a fixed $g\in\FFF$, the branching divisor $B_g$ is supposed to
have a multiple component, say $B_g=rC+R$ with $r>1$.  Let $X'\to Y$
be the standard cyclic covering obtained from $X$ after the
normalization procedure has been applied to the multiple component
$rC$.  If $X'$ is associated to 
\[
  nL'_\chi\equiv\sum_{f\in\FFF}\,[G:\Im\fhat]\, f(\chi)^\bullet B'_f,
\]
then for every $k=1,\ldots,n-1$,
\[
  L'_{\chi^k}\equiv kL_\chi
  -\floor{\frac{krg(\chi)^\bullet}{m_g}}C
  -\floor{\frac{kg(\chi)^\bullet}{m_g}}R
  -\sum_{f\neq g}\floor{\frac{kf(\chi)^\bullet}{m_f}}B_f.
\]
\end{pro}
\proof
If $r\geq m_g$, then $r=qm_g+r_1$, with $0\leq r_1<r$.  The
covering data are modified to 
\begin{equation}\label{eq:data1N}
  L_\chi^{(1)}\equiv L_\chi-qg(\chi)^\bullet C, \quad
  B^{(1)}_g\equiv r_1C+R \quad\text{and}\quad
  B^{(1)}_f\equiv B_f \text{ for } f\neq g. 
\end{equation}
If $(r_1,m_g)=d>1$, then the map
$g_2:\Ghat\arrow{g}\ZZ/m_g\to\ZZ/\frac{m_g}{d}$ is considered.  
The integer $g(\chi)^\bullet$ satisfies 
\begin{equation}\label{eq:rel2N}
  g(\chi)^\bullet = q_1\frac{m_g}{d}+g_2(\chi)^\bullet.
\end{equation}
The covering data are modified to 
\begin{equation}\label{eq:data2N}
  L_\chi^{(2)}\equiv L_\chi^{(1)}-q_1\frac{r_1}{d}C, \quad
  B^{(2)}_g\equiv R, \quad
  B^{(2)}_{g_2}\equiv B^{(1)}_{g_2}+\frac{r_1}{d}C
  \quad\text{and}\quad
  B^{(2)}_f\equiv B^{(1)}_f
  \quad\text{for }
  f\neq g,g_2.
\end{equation}
Finally, if the multiplicity of $C$, $r_1/d$, is an integer
greater than $1$ and prime to $m_g/d$, then the map 
$g_3:\Ghat\arrow{g_2}\ZZ/m_g
  \arrow{r_1/d}\ZZ/\frac{m_g}{d}$ is considered.  
We have 
\begin{equation}\label{eq:rel3N}
  \frac{r_1}{d}g_2(\chi)^\bullet =
  q_2\frac{m_g}{d}+g_3(\chi)^\bullet.
\end{equation}
The covering data are modified to 
\begin{equation}\label{eq:data3N}
  L'_\chi\equiv L_\chi^{(2)}-q_2C, \quad
  B'_{g_2}\equiv B^{(1)}_{g_2}, \quad
  B'_{g_3}\equiv B^{(2)}_{g_3}+C
  \quad\text{and}\quad
  B'_f\equiv B^{(2)}_f
  \quad\text{for }
  f\neq g_2,g_3.
\end{equation}
Using (\ref{eq:data1N}), (\ref{eq:data2N}) and (\ref{eq:data3N}) we
have 
$L'_\chi\equiv L_\chi-
  (qg(\chi)^\bullet+q_1r_1/d+q_2)C$ and, since we
know that 
$L'_{\chi^k}\equiv kL'_\chi-
  \sum\floor{kf(\chi)^\bullet/m_f}B'_f$, 
we also have
\begin{align*}
  L'_{\chi^k}
  &\equiv kL'_\chi 
  -\floor{\frac{kg_2(\chi)^\bullet}{m_g/d}}B_{g_2}
  -\floor{\frac{kg_3(\chi)^\bullet}{m_g/d}}(C+B_{g_3})
  -\floor{\frac{kg(\chi)^\bullet}{m_g}}R
  -\sum_{f\neq g,g_2,g_3}\floor{\frac{kf(\chi)^\bullet}{m_f}}B_f\\
  &\equiv kL_\chi
  -\left(\floor{\frac{kg_3(\chi)^\bullet}{m_g/d}}
  +kqg(\chi)^\bullet+kq_1\frac{r_1}{d}+kq_2
  \right)C
  -\floor{\frac{kg(\chi)^\bullet}{m_g}}R 
  -\sum_{f\neq g}\floor{\frac{kf(\chi)^\bullet}{m_f}}B_f.
\end{align*}
Now, from (\ref{eq:rel3N}) and (\ref{eq:rel2N}), we get successively 
\[
  \floor{\frac{kg_3(\chi)^\bullet}{m_g/d}}
  = \floor{\frac{kr_1g_2(\chi)^\bullet}{m_g}}-kq_2
  = \floor{\frac{kr_1g(\chi)^\bullet}{m_g}}-kq_1\frac{r_1}{d}-kq_2,
\]
and finally, by the Euclidean division of $r$ to $m_g$, 
\[
  \floor{\frac{kg_3(\chi)^\bullet}{m_g/d}}
  = \floor{\frac{krg(\chi)^\bullet}{m_g}}-
  kqg(\chi)^\bullet-kq_1\frac{r_1}{d}-kq_2,
\]
\qed

\begin{pro}
Let $X$ be a normal projective variety and $Y$ be a smooth projective
variety. Let $\pi:X\to Y$ be a cyclic covering.
If $\omega_X$ is a dualizing sheaf for $X$, then
\[
 \pi_\ast\omega_X=\bigoplus_{\chi\in\Ghat}\omega_Y\otimes\Ll_\chi,
\]
the action of $G$ on $\omega_Y\otimes\Ll_\chi$ being the
multiplication by $\chi\inv$. 
\end{pro}
\proof
We recall the following construction from \cite{Ha}, III, Ex.\/6.10
and Ex.\/7.2 valid for $X$ and $Y$ be projective schemes and $\pi:X\to
Y$ a finite morphism.  For any quasi-coherent
$\Oo_Y$-module $\Gg$, the sheaf $\Homcal(\pi_\ast\Oo_X,\Gg)$ is a
quasi-coherent $\pi_\ast\Oo_X$-module.  Hence there exists a unique
quasi-coherent $\Oo_X$-module, denoted $\pi^!\Gg$, such that
$\pi_\ast\pi^!\Gg=\Homcal(\pi_\ast\Oo_X,\Gg)$.  If $\Ff$ is coherent
on $X$ and $\Gg$ is quasi-coherent on $Y$, then there is a natural
isomorphism 
$\pi_\ast\Homcal(\Ff,\pi^!\Gg)\simeq\Homcal(\pi_\ast\Ff,\Gg)$.  
It yields the natural isomorphism
\[
  \Hom(\Ff,\pi^!\Gg) \arrow{\simeq} \Hom(\pi_\ast\Ff,\Gg)
\]
since $H^0(X,\Homcal(\Ff,\pi^!\Gg))\simeq 
  H^0(Y,\pi_\ast\Homcal(\Ff,\pi^!\Gg))$.
If $\omega_Y$ is the canonical sheaf for $Y$, then it follows that
$\pi^!\omega_Y$ is a dualizing sheaf for $X$.  Hence 
\[
  \pi_\ast\omega_X=\pi_\ast\pi^!\omega_Y
  =\Homcal(\pi_\ast\Oo_X,\omega_Y)
  =\bigoplus_{\chi\in\Ghat}\omega_Y\otimes\Ll_\chi.
\]
\qed

\begin{lem}\label{l:formulaForq}
Let $S_1\to Y$ be a normal standard cyclic covering of surfaces
defined by the line bundle $\Ll$. 
If $S_1$ has only rational singularities and $S\to S_1$ denotes a 
desingularization of $S_1$, then 
\[
  q(S) = q(Y)+\sum_{k=1}^{n-1}h^1(Y, \omega_Y\otimes\Ll_{\chi^k})
       -\sum_{k=1}^{n-1}h^2(Y, \omega_Y\otimes\Ll_{\chi^k}).
\]
\end{lem} 
\proof
Since the singularities are rational, if 
$S\arrow{\epsilon}S_1$ is a resolution of the singular points 
of $S_1$, then $R^i\epsilon_\ast\Oo_S=0$, for all $i\geq1$. 
From the Leray spectral sequence it follows that
$h^i(S,\Oo_S)=h^i(S_1,\Oo_{S_1})$ for all $i$, and hence
$\chi(\Oo_S)=\chi(\Oo_{S_1})$.  Then
$q(S)=q(S_1)=p_g(S_1)+1-\chi(\Oo_{S_1})=
  h^0(Y,\pi_\ast\omega_{S_1})+1-\chi(\pi_\ast\Oo_{S_1})$ and 
using the formulae for $\pi_\ast\omega_S$ and $\pi_\ast\Oo_S$, we get 
\[
  q(\Stilde)= \sum_{k=0}^{n-1}h^0(Y,\omega_Y\otimes\Ll_{\chi^k})+1-
        \sum_{k=0}^{n-1}\chi(\Ll_{\chi^j}^{-1}).
\]
By Serre duality, the required equality follows.
\qed

One more notation is in order.  Let $P$ be a singular point of $B$ and
let $\mu:Y\to\PP^2$ be a desingularization of $B$ at $P$, with
$E_{P,1},E_{P,2},\ldots$\/ be the irreducible components of the fibre
$\mu\inv(P)\subset Y$.  $\eE_{P}$ will denote this finite array of
irreducible curves, and if $\cc$ is a finite array of rational numbers
$c_1,c_2,\ldots$, then
\begin{equation}
  \label{eq:cE}
  \cc\cdot\eE_P = \sum_\alpha c_\alpha E_{P,\alpha}.
\end{equation}

\likeproof[Proof of Theorem \ref{th:q}] 
For any integer $n$, $S_0\subset\PP^3$, the $n$-cyclic multiple plane
associated to $B$ and $H_\infty$ is considered.  $B$ is assumed to be
reduced and transverse to $H_\infty$.  By Proposition
\ref{p:multipleToCovering} there is a convenient sequence of
blowing-ups such that $S_1$, the strict transform of $S_0$, becomes a
standard cyclic covering of the plane defined by 
$nL'_\chi\equiv B+(\beta n-b)H_\infty$, with $\beta=\ceil{b/n}$ and
$L'_\chi\equiv \beta H$.  We choose a desingularization of $B$ such
that its total transform by $\mu:Y\to\PP^2$ is a divisor with normal
crossing intersections.  $S_1$ induces a standard cyclic covering
$S_2$ defined by
\[
  nL''_\chi \equiv B+(\beta n-b)H_\infty+\sum_P\cc_P\cdot\eE_P.
\]
We apply the normalization procedure to $S_2$ to end up with a
normalization $S$ of $S_0$ that has only Hirzebruch-Jung singularities
(see \cite{Pa}, Proposition 3.3).  By Proposition
\ref{p:keyProposition}, if $S$ is defined by the line bundle
$\Ll_\chi$, then 
\begin{align}
  L_{\chi^k}
  &\equiv kL''_\chi-\floor{\frac{k}{n}(\beta n-b)}H-
  \sum_P\floor{\frac{k}{n}\cc_P}\cdot\eE_P \notag\\
  &\equiv
  \ceil{\frac{kb}{n}}H-\sum_P\floor{\frac{k}{n}\cc_P}\cdot\eE_P,
\label{eq:Lchik}
\end{align}
the last equality resulting from
$\beta k-\floor{k(\beta n-b)/n}=\ceil{kb/n}$.  Here,
$\floor{k\cc_P/n}\cdot\eE_P$ denotes 
$\sum_\alpha\floor{kc_{P,\alpha}/n}E_{P,\alpha}$.
From Lemma~\ref{l:formulaForq}, since
$H\cdot(-L_{\chi^k})=-\ceil{kb/n}<0$, it follows that 
\begin{equation}\label{eq:q}
  q(S) = \sum_{k=1}^{n-1}h^1(Y,K_Y+L_{\chi^k}). 
\end{equation}

In order to end the proof we have to take account in the formula
above, of the vanishing of certain $h^1$'s and of the equality of the
others with certain superabundances of linear systems on the
projective plane.

\paragraph{Claim}
$H^1(Y,\omega_Y\otimes\Ll_{\chi^k})
  \simeq H^1(\PP^2,\Oo_{\PP^2}(-3+\ceil{kb/n})
  \otimes\Ii_{Z(k/nB)})$,
with $Z(k/nB)$ the scheme defined by the multiplier
ideal of $k/n B$.

Indeed, by (\ref{eq:Lchik}) and (\ref{eq:upDown}), it follows
that
\begin{align*}
  H^1(Y,\omega_Y\otimes\Ll_{\chi^k}) 
  &= H^1(Y,\mu^\ast\omega_{\PP^2}\otimes\Oo_Y(\ceil{\frac{kb}{n}}H)
  \otimes\Oo_Y(K_{Y|\PP^2}-\sum_P\floor{\frac{k}{n}\cc_P}\cdot\eE_P))\\
  &= H^1(Y,\mu^\ast\omega_{\PP^2}\otimes\Oo_Y(\ceil{\frac{kb}{n}}H)
  \otimes\Oo_Y(K_{Y|\PP^2}-\floor{\frac{k}{n}}B
  -\sum_P\floor{\frac{k}{n}\cc_P}\cdot\eE_P))\\
  &= H^1(Y,\mu^\ast\omega_{\PP^2}\otimes\Oo_Y(\ceil{\frac{kb}{n}}H)
  \otimes\Oo_Y(K_{Y|\PP^2}-\floor{\mu^\ast\frac{k}{n}B}))\\
  &\simeq H^1(\PP^2,\Oo_{\PP^2}(-3+\ceil{\frac{kb}{n}})
  \otimes\Ii_{Z(\frac{k}{n}B)}),
\end{align*}
justifying the claim.

\medskip

Using (\ref{eq:q}) and the above claim, the irregularity is given by 
\begin{equation*}
  \label{eq:qInGeneral}
  q(S) 
  = \sum_{k=1}^{n-1}h^1(\PP^2,\Oo_{\PP^2}(-3+\ceil{\frac{kb}{n}})
  \otimes\Ii_{Z(\frac{k}{n}B)}).
\end{equation*}
If $k/n\not\in J(B,n)$, then either $k/n$ is not a jumping number of
$B$, or it is, but $kb/n$ is not an integer.  In the former case, if
$\xi$ is the biggest jumping number for $B$ smaller than $k/n$, then,
since $\ceil{kb/n}-\xi>0$,
\[
  h^1(\PP^2,\Oo_{\PP^2}(-3+\ceil{\frac{kb}{n}})
  \otimes\Ii_{Z(\frac{k}{n}B)}) 
  = h^1(\PP^2,\Oo_{\PP^2}(-3+\ceil{\frac{kb}{n}})
  \otimes\Ii_{Z(\xi B)}) 
  = 0
\]
by Kawamata-Viehweg-Nadel Vanishing Theorem.  In the latter case, we
apply the same argument, now using $\ceil{kb/n}-kb/n>0$.  The result
follows. 
\qed

\begin{cor} \label{c:qIrr}
Under the hypotheses of Theorem \ref{th:q}, if furthermore $B$ is
supposed to be an irreducible plane curve, then
\[
  q(S)
  = \sum_{\xi\in J'(B,n)}
  h^1(\PP^2,\Ii_{Z(\xi B)}(-3+\xi b)),
\]
with $J'(B,n)$ the subset of $J(B,n)$ that contains those rationals
$\xi$ for which the denominator can not be the power of a prime. 
\end{cor}
\proof
In \cite{Za2} the following topological result is established: 
{\it If $q$ is the power of a prime and $B$ is
irreducible and transverse to $H_\infty$, then the $q$-cyclic multiple
plane is regular.} 
By inspecting the formula for the irregularity given in Theorem
\ref{th:q} for $q$-multiple planes associated to $B$ and $H_\infty$,
$q$ a power of a prime such that there exists a jumping number 
$l/q\in J(B,n)$, we obtain the corollary.
\qed

\section{The case of specified singularities}
  \label{s:clusters} 
Explicit versions of Theorem \ref{th:q} may be formulated as soon as
the multiplier ideals and the jumping numbers can be evaluated.  Such
an explicit version is obtained, for example, if the singularities of
$B$ are locally characterized by the equation $x^{dp}-y^{dq}=0$, with
$p,q$ and $d$ positive integers and $p,q$ relatively prime.  Turning
to Definition \ref{d:minCluster} where the cluster
$K_{p,q}(\alpha,\beta)$ is introduced and if $Z_{p,q}(\alpha,\beta)$
is the subscheme associated to it, we have:

\begin{cor} \label{c:pqTypeSing}
Let $B$ be a plane curve of degree $b$ with each of its singular
points of type, either $A_1$, or given locally by the equation
$x^{dp}-y^{dq}=0$.  Let $H_\infty$ be a line transverse to $B$ and let
$S$ be a desingularization of the $n$-cyclic multiple plane associated
to $B$ and $H_\infty$.  Then
\[
  q(S) 
  = \sum_{(\alpha,\beta)}h^1(\PP^2,
  \Ii_{\Zz_{p,q}(\widetilde{\alpha,\beta})}(-3+\frac{\alpha p+\beta q}{dpq}b)).
\]
The sum ranges over the couples $(\alpha,\beta)$ such that
$\dfrac{\alpha p+\beta q}{dpq}<1$ and 
$\dfrac{\alpha p+\beta q}{dpq}\in\dfrac{1}{\gcd(b,n)}\ZZ$.  In
addition, the couple of positive integers $(\widetilde{\alpha,\beta})$
is defined by
\[
  \min_{(\alpha',\beta')}
  \{\alpha'p+\beta'q\geq (\alpha-1)p+(\beta-1)q+1\},
\]
and 
$\Zz_{p,q}(\widetilde{\alpha,\beta})= \cup_P
Z_{p,q}(\widetilde{\alpha,\beta})_P$, 
with $P\in\Sing C$ not of type $A_1$.
\end{cor}

In this section we mainly want to establish this corollary.  
We need to control the jumping numbers and the multiplier ideals
associated to a curve with this type of singularities.  
The multiplier ideals and their jumping numbers are known in this
case; see for example \cite{Ei} and \cite{EiLa}, or \cite{Ho} for the
case of monomial ideals in general.  We like to present a different
argument based on Enriques diagrams for the particular case of two
unknowns, since it will provide a simple interpretation of the
multiplier ideals involved, and also, could provide an algorithm for
the generalization to an arbitrary singular point of a curve on a
surface.

\subsection{Clusters and Enriques diagrams} 
Let $X$ be a surface and $P\in X$ a smooth point.  A point $Q$ is
called infinitely near to $P$ if $Q\in X'$, $\mu:X'\to X$ is a
composition of blowing-ups and $Q$ lies on the exceptional
configuration that maps to $P$.

\begin{definition*}
A {\it cluster} in $X$, centered at a smooth point $P$ is a finite set
of weighted infinitely near points to $P$,
$K=\{P_1^{w_1},\ldots,P_r^{w_r}\}$, with $P_1=P$.
\end{definition*}

Let $\mu:Y\to X$ be the composition of blowing-ups 
$Y=Y_{r+1}\to Y_r\to\cdots\to Y_1=X$, with
$Y_{\alpha+1}=\Bl_{P_\alpha}Y_\alpha$.  Since the points infinitely
near $P$ are partially ordered---the point $Q$ precedes the point $R$
if and only if $R$ is infinitely near $Q$\/---the points of a cluster
are partially ordered.  In the sequel, if $K$ is a cluster, then all
points preceding a point that belongs to $K$ are in $K$, possibly with
weight $0$.

Let $K$ be a cluster centered at $P$.  Each point $P_\alpha$
corresponds to an exceptional divisor $E_\alpha\subset Y_{\alpha+1}$.  
All its strict transforms will also be denoted by $E_\alpha$ and the
total transform of each $E_\alpha$ will be denoted by $F_\alpha$.
When needed, the strict transform of $E_\alpha$ on $Y_\beta$ will be
denoted by $E_\alpha^{(\beta)}$, and similarly for the total
transform.  For example $F_\alpha^{(\alpha+1)}=E_\alpha^{(\alpha+1)}$.

Every cluster $K$ defines a divisor $D_K=\sum w_\alpha F_\alpha$ on
$Y$ and an ideal sheaf $\mu_\ast\Oo_Y(-D_K)$ on $X$, hence a subscheme
$Z_K$ of $X$.  The lemma below clarifies the comparison between the
ideal sheaf $\Oo_Y(-D_K)$ and the pull-back
$\mu^\ast\mu_\ast\Oo_Y(-D_K)$.

\begin{definitions*} 
Let $K$ be a cluster.
A point $P_\beta$ is said to be {\it proximate} to $P_\alpha$ if 
$P_\beta$ lies on $E_\alpha$, the exceptional divisor corresponding
to the blowing-up at $P_\alpha$, or on one of its strict transforms. 

A cluster $K$ is said to satisfy the proximity relations if for every
$P_\alpha$ in $K$,
\[
  \wbar_\alpha = 
  \sum_{P_\beta \text{ proximate to }P_\alpha}w_\beta \leq w_\alpha.
\]
\end{definitions*}

\begin{lem}[see also \cite{Ca}, Theorem 4.2]
  \label{l:unloadingProcedure} 
Let $K=\{P_1^{w_1},\ldots,P_r^{w_r}\}$ be a cluster that
contains a point $P_\alpha$ for which the proximity relation is not
satisfied.  If $K'=\{P_1^{w'_1},\ldots,P_r^{w'_r}\}$ is the
cluster defined by $w'_\alpha=w_\alpha+1$, $w'_\beta=w_\beta-1$ for
every $\beta$ with $P_\beta$ proximate to $P_\alpha$, and
$w'_\gamma=w_\gamma$ otherwise, then $K$ and $K'$ define the same
subscheme in $X$, \ie $\mu_\ast\Oo_Y(-D_K)=\mu_\ast\Oo_Y(-D_{K'})$. 
\end{lem}
$K'$ is said to be obtained from $K$ by the unloading procedure.
Starting from $K$, iterated applications of this procedure lead to a
cluster $\Ktilde$ that satisfies the proximity relations and defines
the same subscheme in $X$.  $\Ktilde$ is called the unloaded
cluster.  Notice that 
\[
  \mu^\ast\mu_\ast\Oo_Y(-D_K)
  \simeq \mu^\ast\mu_\ast\Oo_Y(-D_\Ktilde) 
  \simeq \Oo_Y(-D_{\Ktilde}).
\]

\begin{rem*}
If $w_r<0$, then the proximity relation is not satisfied
at $P_r$ since $\wbar_r=0$.
When the unloading procedure of Lemma \ref{l:unloadingProcedure} is
applied to a cluster with non-negative weights, it may happen that a
weight becomes negative, or more precisely, becomes $-1$.  But it is
to be noticed that the negative weight is eventually rubbed out by the
next applications of the procedure, and that the unloaded cluster
has only non-negative weights.  Moreover, the unloaded cluster
associated to a cluster with non-positive weights is the empty
cluster, the one with all its weights equal to $0$.
\end{rem*}

\begin{definition*}
A {\it gridded tree} is a couple $(T,g)$, where $T=T(\VVV,\AAA)$ is an
oriented tree with $\VVV$ the set of vertices and $\AAA$ the set of
arcs, and $g$ is a map 
\[
g:\AAA\to\{\text{slant}, \text{horizontal},\text{vertical}\}.
\] 
\end{definition*}

\begin{definition*}
Let $T$ be a gridded tree.  A horizontally (vertically) 
{\it $L$-shape branch}
of $T$ is an ordered chain of arcs, such that each begins where the
previous ends, and such that all are horizontal (vertical), but the
first.   
\end{definition*}

Notice that an $L$-shape branch is completely determined by the subset
of incident vertices of it.  Moreover, an arc is an $L$-shape branch,
regardless its value through $g$.

\begin{definition*}
Let $T$ be a gridded tree.  A {\it segment} is a maximal chain of arcs
of the same type through $g$, arcs that are also maximal $L$-shape
branches. 
\end{definition*}

\begin{example} \label{ex:tpq}
Let $p<q$ be relatively prime positive integers.  $T_{p,q}$ will
denote the gridded tree associated to the Euclidean algorithm.  If 
$r_0=a_1r_1+r_2,\ldots,r_{m-2}=a_{m-1}r_{m-1}+r_{m}$ and
$r_{m-1}=a_mr_m$, with $r_0=q$ and $r_1=p$, then $T_{p,q}$ has $d$
segments containing $a_1,\ldots,a_{m-1}$ and respectively $a_m$
vertices each.  The first segment is slanted and the others are
alternatively, either horizontal or vertical, starting with a
horizontal one. 
\end{example}

\begin{definition*}
An {\it Enriques diagram} is an weighted gridded tree. 
\end{definition*}

Clusters and Enriques diagrams carry the same information as the lemma
below asserts, and it will often be convenient to argue using
diagrams.

\begin{definition*}
A point of a cluster is said to be {\it free} if it is proximate to
exactly one point of the cluster.  A point is said to be a 
{\it satellite} if it is proximate to exactly two points of the
cluster. 
\end{definition*}

\begin{lem}[see \cite{Ev}]
There exists an unique map from the set of clusters in $X$ centered at
a smooth point $P$ to the set of Enriques diagrams such that:
\begin{enumerate}\setlength{\itemsep}{-.25\baselineskip}
\item 
for every cluster $K=\{P_1^{w_1},\ldots,P_r^{w_r}\}$ the set of
vertices of the image tree is $\VVV=\{P_1,\ldots,P_r\}$ with the
weights given by the integers $w_1,w_2,\ldots,w_r$; 
\item
at every point ends at most one arc;
\item
a point $P_\alpha$ is satellite if and only if there is either a
horizontal or a vertical arc that ends at the vertex $P_\alpha$;
\item
if there is an arc that begins at the vertex $P_\alpha$ and ends at
the vertex $P_\beta$ then $P_\beta\in E_\alpha^{(\beta)}$, and the
converse is true if $P_\beta$ is free;
\item
$P_\beta$ is proximate to $P_\alpha$ if and only if there is an
$L$-shape branch that starts at $P_\alpha$ and ends at $P_\beta$;
\item
the strict transforms $E_\alpha$ and $E_\beta$ intersect on $Y$ if and
only if the Enriques diagram contains a maximal $L$-shape branch that
has $P_\alpha$ and $P_\beta$ as its extremities; 
\item 
an arc that begins at a vertex of a free point and ends at a vertex of
a satellite point is horizontal. 
\end{enumerate}
\end{lem}

\subsection{The minimal unloaded clusters associated to a $T_{p,q}$
  tree}

Let $p<q$ be relatively prime positive integers.  All clusters treated
in this subsection will be associated to the gridded tree $T_{p,q}$
introduced in Example \ref{ex:tpq}.  The intent is to look for a
characterization of the minimal unloaded clusters modeled on
$T_{p,q}$.  We refer to Lemma \ref{l:minCluster} for the result. 
 
Depending on the context its vertices will be denoted either by
$P_\alpha$ \ie using one
subscript $1\leq\alpha\leq r=a_1+\cdots+a_m$, or by $P_{k,i}$, \ie
using two subscripts $1\leq k\leq d$, $1\leq i\leq a_k$.
$T_{3,5}$ is represented below in the latter notation.

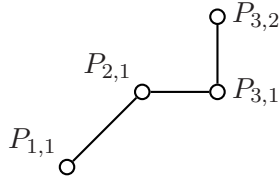
\begin{center}
\begin{pspicture}(-2,-1)(2,1.25)%\showgrid
  \psset{radius=.11}
  \Cnode(-1,-1){1}
  \Cnode(0,0){2}
  \Cnode(1,0){3}
  \Cnode(1,1){4}
  \ncline{1}{2}
  \ncline{2}{3}
  \ncline{3}{4}
  \uput{\labelsep}[ul](-1,-1){$P_{1,1}$}
  \uput{\labelsep}[ul](0,0){$P_{2,1}$}
  \uput{\labelsep}[r](1,0){$P_{3,1}$}
  \uput{\labelsep}[r](1,1){$P_{3,2}$}
\end{pspicture} 
\end{center}

Let $K$ be a cluster.  We define the proximity matrix of $K$ by 
$\Pi=||p_{\alpha\beta}||$, where the elements of the diagonal equal
$1$ and, for every $\alpha\neq\beta$, the element $p_{\alpha\beta}$
equals $-1$ if $P_\beta$ is proximate to $P_\alpha$ and $0$ if not.  
Notice that along the $\alpha$ column of $\Pi$, the non-zero elements
not on the diagonal correspond to the points to which $P_\alpha$ is a
satellite.  

The proximity matrix is the decomposition matrix of the strict
transforms $E_\alpha$'s in terms of the total transforms
$F_\alpha$'s.  Hence if $K=\{P_1^{w_1},\ldots,P_r^{w_r}\}$, then on
$Y$,
\[
  D_K = \sum_\alpha w_\alpha F_\alpha 
  = \sum_\alpha c_\alpha E_\alpha,
\]
and $\cc=\ww\,\Pi\inv$, where $\ww=(w_1,\ldots,w_r)$ and similarly
$\cc=(c_1,\ldots,c_r)$.  The formula
\[
  E_\alpha 
  = F_\alpha-\sum_{P_\beta\text{ proximate to }P_\alpha}F_\beta
\]
and induction on $\alpha$ tell us that the coefficient of $E_r$ in the
decomposition of a total transform corresponding to a point lying on
the $k$th segment in terms of strict transforms equals the remainder
$r_k$ introduced in Example \ref{ex:tpq}.

\begin{lem}\label{l:ap+bq}
If $K=\{P_1^{w_1},\ldots,P_r^{w_r}\}$ is an unloaded cluster centered
at $P$, then the coefficient of $E_r$ is of the form $ap+bq$, with
$a,b$ non-negative integers. 
\end{lem}
\proof
We shall denote the weights on the $k$th segment of the Enriques
diagram for the cluster $K$ by $w_{k,1},w_{k,2},\ldots,w_{k,{a_k}}$
and the coefficient of $E_r$ by $c_r$.  
We shall successively transform the cluster, each time considering the
last segment that contains non-zero weights, unless this segment is
the first or the second one.  The transformation is the following: if
the last segment with non-zero weights is the segment $k+1$, and if 
$\wbar_{k,a_k}=\sum_iw_{k+1,i}$, then put 
\begin{enumerate}\setlength{\itemsep}{0pt}
\item
  $w'_{k+1,i}=0$ for $1\leq i\leq a_{k+1}$,
\item 
  $w'_{k,i}=w_{k,i}-\wbar_{k,a_k}$ for $1\leq i\leq a_k$,
\item
  $w'_{k-1,1}=w_{k-1,1}+\wbar_{k,a_k}$,
\end{enumerate}
and leave the other weights unchanged.  It is easy to see that the
cluster $K'$ defined by the weights $w'_{k,i}$ is again unloaded
and that the coefficient of $E_r$ remains unchanged.  
Hence the same process can be applied till eventually $K$ is
transformed into the cluster $K'$ with non-zero weights only on the
first and, at the most, the second segments of the Enriques diagram.
$K'$ is unloaded and 
\[
  c_r = c'_r 
  = p\sum_{i=1}^{a_1}w'_{1,i} + r_2\sum_{j=1}^{a_2}w'_{2,j}
  = p\sum_{i=1}^{a_1}(w'_{1,i}-
  \overline{w'}_{1,a_1}) + q\overline{w'}_{1,a_1},
\]
where $q=a_1p+r_2$.
\qed

Besides, since for each couple of non-negative integers $a,b$ there exists
an unloaded cluster with $ap+bq$ the coefficient of $E_r$---for example
the cluster with $w_{1,1}=a+b$, $w_{1,i}=b$ for $i\neq1$, $w_{2,1}=b$
and $w_{k,i}=0$ otherwise---, 
the ideal sheaf $\mu_\ast\Oo_Y(-(ap+bq)E_r)$ defines the minimal
unloaded cluster having $ap+bq$ the coefficient of $E_r$.  It is
natural to ask the question whether we can decide if an unloaded
cluster is minimal only by inspection of its weights, or equivalently
its associated divisor.  The answer is yes and is given by the lemma
hereafter.  It will deal with clusters satisfying the following
condition:
\begin{verse}
 ($\ast$) 
for every ordered chain of maximal $L$-shape branches determined by
the points $P_{\alpha_1},\ldots,P_{\alpha_{l}}$, \ie each
$P_{\alpha_j}$ precedes $P_{\alpha_{j+1}}$ and the $j$th maximal
$L$-shape branch starts at $P_{\alpha_j}$ and ends at
$P_{\alpha_{j+1}}$, then  
\[
  \sum_{j=1}^l(w_{\alpha_j}-\wbar_{\alpha_j})
  < \sum_{j=1}^lp_{\alpha_j}+2-l.
\]
\end{verse}

\begin{definition} \label{d:minCluster}
If $a$ and $b$ are non negative integers, $K_{p,q}(a,b)$ denotes the
minimal unloaded cluster associated to the $T_{p,q}$ tree and
whose coefficient of the last strict transform equals $ap+bq$. 
\end{definition}

\begin{lem} \label{l:minCluster}
Let $K$ be an unloaded cluster with $ap+bq$ the coefficient of its
last strict transform.  $K$ satisfies ($\ast$) if and only if
$K=K_{p,q}(a,b)$.
\end{lem}
\proof
We start by showing that a minimal unloaded cluster $K^\text{min}$
always satisfies ($\ast$).  Indeed, if not, there would exist a chain
of maximal $L$-shape branches such that 
$p_{\alpha_j}\geq w_{\alpha_j}-\wbar_{\alpha_j}$ for $1\leq j\leq l$,
and that 
$\sum_1^lw_{\alpha_j}-\wbar_{\alpha_j}\geq\sum_1^l(p_{\alpha_j}-1)+2$.  
Furthermore, since there would be at least two points such that
$p_\alpha=w_\alpha-\wbar_\alpha$, we may always assume 
$w_{\alpha_1}-\wbar_{\alpha_1}=p_{\alpha_1}$, 
$w_{\alpha_j}-\wbar_{\alpha_j}=p_{\alpha_j}-1$ 
for $2\leq j\leq l-1$, and
$w_{\alpha_l}-\wbar_{\alpha_l}=p_{\alpha_l}$.  
Now for this chain, we could apply the inverse of the unloading
procedure successively at
$P_{\alpha_1},P_{\alpha_2},\ldots,P_{\alpha_l}$ to end up with an
unloaded cluster with the same coefficient for $E_r$ as
$K^\text{min}$, hence a contradiction.

To end the proof, we assume that $K\neq K^\text{min}$, $K$ and
$K^\text{min}$ having the same coefficient for $E_r$, and show that
$K$ does not satisfy ($\ast$).
Since it is satisfied by $K^\text{min}$, we notice that along
each segment of $K^\text{min}$, there is at most one jump of height
$1$.   We may further assume that $w_{1,1}\geq w_{1,1}^\text{min}+1$.
To make up for the apparent increase of $c_{m,a_m}\,(=c_r)$ 
by at least $r_1=p$ due to the difference between $w_{1,1}$ and 
$w_{1,1}^\text{min}$, some of the weights along the next segments of
$K$ must be smaller than the corresponding weights of $K^\text{min}$,
but not along the first segment.  Looking at the points 
on the second segment and at the first point of the third segment for
this counterbalance problem, we notice that at most one of their
weights may not diminish, otherwise ($\ast$) will not be satisfied
somewhere along the following segments.  Two possibilities can
appear. 
First, all $a_2$ weights of the second segment satisfy
$w_{2,i}\leq w_{2,i}^\text{min}-1$ and, either there exists an $i$ such
that $w_{3,i}\leq w_{3,i}^\text{min}-1$, or
$w_{3,i}=w_{3,i}^\text{min}$ for all $i$'s.  In the former case
($\ast$) is not verified at
$P_{1,a_1},P_{3,1},P_{3,2},\ldots,P_{3,i}$, and in the latter the same
counterbalance problem must be solved for a difference of $r_3$ units,
starting with the fourth segment.  
Second, the inequalities $w_{2,i}\leq w_{2,i}^\text{min}-1$ are
verified for all but one point of the second segment,
and $w_{3,1}\leq w_{3,1}^\text{min}-1$.  There is a counterbalance
problem left for $r_2$ units, starting with the third segment. 
Eventually, the counterbalance problem is pushed on to the last
segment and hence ($\ast$) will not be satisfied there for $K$.  
\qed

\subsection{The multiplier ideals and the jumping numbers for
  $(x^{dp},y^{dq})$} 

Let $P$ be a singular point of $B\subset X$ given locally by
$x^{dp}+y^{dq}=0$, with $p,q$ and $d$ positive integers and $p\leq q$
relatively prime, and let $\mu:Y\to X$ be the minimal log resolution
of $B$ at $P$.  The exceptional configuration of $\mu$ is given by the
gridded tree $T_{p,q}$.  As before, we shall denote by $E_r$ the last
strict transform.

\begin{lem}\label{l:lastCoeff}
The coefficient of $E_r$ in
$-K_{Y|X}+\floor{\mu^\ast \xi B}$ equals $\floor{dpq\xi}-(p+q-1)$. 
\end{lem}
\proof
It is sufficient to determine the coefficient of $E_r$ in 
$\mu^\ast B$.  By Example \ref{ex:tpq} and the decomposition of the
$F_\alpha$'s in terms of the strict transforms, we have that the
coefficient of $E_r$ is $a_1r_1^2+\cdots+a_mr_m^2=dpq$. 
\qed

\begin{pro}\label{p:unloadingOnLastCoeff}
If $c_r$ is the coefficient of $E_r$ in 
$-K_{Y|X}+\floor{\mu^\ast \xi B}$, then the multiplier ideal 
$\Jj(\xi B)$ is given by 
\[
  \Jj(\xi B) = \mu_\ast\Oo_Y(-c_rE_r),
\]
\ie is the ideal sheaf associated to the minimal cluster that contains
$c_rE_r$. 
\end{pro}
\proof
We shall argue on the cluster associated to the divisor
$-K_{Y|X}+\floor{\mu^\ast \xi B}$.  To find the multiplier ideal is
equivalent to determine the unloaded corresponding cluster.  Let the
pull-back of $B$ be $\sum_1^rc_\alpha E_\alpha + B = \cc\cdot\eE +B$.
Then 
$-K_{Y|X}+\floor{\mu^\ast \xi B} = \sum_1^rw_\alpha F_\alpha 
  = \ww\cdot\fF$,
with $\ww = -\oomega + \floor{\xi\cc}\cdot\Pi$ and
$\oomega=(1,\ldots,1)$. 

Let $P_{\alpha_1},\ldots,P_{\alpha_l}$ be ordered points that
determine a chain of maximal $L$-shape branches. 
Since 
\begin{equation}\label{eq:proximityVector}
  \ww-\overline{\ww} = \ww\,^t\Pi =
  \floor{\xi\cc}\Pi\,^t\Pi-\oomega\,^t\Pi,
\end{equation}
where the matrix $-\Pi\,\Pi^t$ is the intersection matrix of the
strict transforms $E_\alpha$'s on the surface $Y$, for every
$1\leq j\leq l$, 
\[
  w_{\alpha_j}-\wbar_{\alpha_j} 
  = -\floor{\xi c_{\alpha_{j-1}}}+
    (p_{\alpha_j}+1)\floor{\xi c_{\alpha_j}}-
    \floor{\xi c_{\alpha_{j+1}}}+(p_{\alpha_j}-1).
\]
So $\sum_{j=1}^l(w_{\alpha_j}-\wbar_{\alpha_j})$ equals
\[
  -\floor{\xi c_{\alpha_0}}
  +p_{\alpha_1}\floor{\xi c_{\alpha_1}}
  +\sum_{j=2}^{l-1}(p_{\alpha_i}-1)\floor{\xi c_{\alpha_j}}
  +p_{\alpha_l}\floor{\xi c_{\alpha_l}}
  -\floor{\xi c_{\alpha_{l+1}}}
  +\sum_{j=1}^l(p_{\alpha_j}-1),
\]
and since $\cc\,\Pi\,^t\Pi=(0,\ldots,0,d)$, we have
\begin{equation}
  \label{eq:ast}
  -2 < \sum_{j=1}^l(w_{\alpha_j}-\wbar_{\alpha_j})
  < \sum_{j=1}^lp_{\alpha_j}+2-l.
\end{equation}
Putting $l=1$ we observe that if the proximity relation is not
satisfied at $P_\alpha$, then $w_\alpha-\wbar_\alpha=-1$.   
But the unloading procedure of Lemma \ref{l:unloadingProcedure} at
$P_\alpha$ changes the vector $\ww-\overline{\ww}$ into the vector 
$\ww-\overline{\ww}+(\Pi\,^t\Pi)_\alpha$.  It follows that the
unloading procedure does not change the inequalities in (\ref{eq:ast})
for the new cluster.  We conclude that the unloaded cluster satisfies
($\ast$), hence, by Lemma \ref{l:minCluster}, the result.
\qed

\begin{pro} \label{p:jumpingNumbers}
The jumping numbers of $B$ at $P$ are $(ap+bq)/(dpq)$ with $a,b$
positive integers. 
\end{pro}
\proof
Let $\xi$ be a jumping number and $K'$ the corresponding unloaded
cluster and let, by Lemma \ref{l:ap+bq}, $c'_r=a'p+b'q$ be its
coefficient for $E_r$ with $a',b'$ non-negative integers.  
By Lemma \ref{l:lastCoeff} and Proposition
\ref{p:unloadingOnLastCoeff}, we have 
\[
  \floor{dpq\xi}-(p+q-1) = ap+bq+1 \leq a'p+b'q
\]
where $a$ and $b$ are non-negative integers, and $a'p+b'q$ is the
first integer combination of $p$ and $q$ with this property.  So, by
the definition of the jumping numbers, $\xi=((a+1)p+(b+1)q)/(dpq)$. 
\qed

\likeproof[Proof of Corollary \ref{c:pqTypeSing}]
We know that the irregularity is given by 
\[
  \sum_{\xi\in J(B,n)}h^1(\PP^2,\Ii_{Z(\xi B)}(-3+\xi b))
\]
with $J(B,n)$ the subset of jumping numbers $\xi$ of $B$ of the form 
$k/n$, $0<k<n$, and such that $\xi b$ is an integer.  By Proposition
\ref{p:jumpingNumbers}, it is sufficient to describe the subscheme
associated to the multiplier ideal $\Jj(\xi B)$ for every 
$\xi=(\alpha p+\beta q)/(dpq)\in J(B,n)$.
%Clearly such a subscheme is supported at the singular points not of
%type $A_1$.  
By Proposition \ref{p:unloadingOnLastCoeff}, Lemma \ref{l:lastCoeff}
and Lemma \ref{l:minCluster}, the subscheme is given by the minimal
unloaded cluster whose coefficient for the last strict transform is
the first integer combination of $p$ and $q$ not smaller than
$(\alpha-1)p+(\beta-1)q+1$.  \qed

\begin{remark} \label{r:oddEven}
Since many of the applications in the next section will be for a curve
$B$ with singularities of a given type $A_\bullet$, we interpret
Corollary \ref{c:pqTypeSing} for this situation.  If $P_1,\ldots,P_r$
are the infinitely near points to $P=P_1$ involved in the minimal log
resolution of an $A_\bullet$ type singularity at $P$, we shall denote by
$Z_P^{[\alpha]}$ the curvilinear subscheme associated to the unloaded
cluster $\{P_1,\ldots,P_\alpha\}$, and by 
$\Zz^{[\alpha]}=\cup Z_P^{[\alpha]}$.

i) If the singularities of $B$ are of type $A_1$ or $A_{2r-1}$,
\ie $p=1$, $q=r$ and $d=2$, then
\[
  q(S) =
  \sum_\alpha
  h^1(\PP^2,\Ii_{\Zz^{[\alpha]}}(-3+\frac{(\alpha+r)b}{2r})),
\]
$\alpha$ ranging from $1$ to $r-1$ such that
$\frac{\alpha+r}{2r}\in\frac{1}{\gcd(b,n)}\,\ZZ$. 

ii) If the singularities of $B$ are of type $A_1$ or $A_{2m}$,
\ie $p=2$, $q=2m+1$ and $d=1$, then
\[
  q(S) =
  \sum_\alpha
  h^1(\PP^2,\Ii_{\Zz^{[\alpha]}}
  (-3+\frac{\alpha b}{2m+1}+\frac{b}{2})),
\]
$n$ and $b$ are even, and $\alpha$ ranges from $1$ to $m$ such that
$\frac{\alpha}{2m+1}\in\frac{1}{\gcd(b,n)}\,\ZZ$.
\end{remark}

\section{Applications}
  \label{s:applications}
We shall now apply the results in the previous sections to
illustrate how to compute in an uniform way, the irregularity for some
examples of cyclic multiple planes.

\subsection*{Zariski's example}
The curve $B$ is irreducible, of degree $6$ and has six cusps as
singularities.  In the formula for the irregularity of the $6$-cyclic
multiple plane in Remark \ref{r:oddEven} ii), since $m=1$, $\alpha$
may only be $1$.  Hence $q(S)=h^1(\PP^2,\Ii_{\Zz}(2))$, where $\Zz$ is
the support of the cusps.  So either the cusps lie on a conic and the
irregularity is $1$, or they do not, and the irregularity is
$0$.  Notice that the result is the same for every $n$-cyclic multiple
plane, provided that $6$ divides $n$.

\subsection*{Artal-Bartolo's first example in \cite{Ar}}
Let $C\subset\PP^2$ be a smooth elliptic curve and let $P_1,P_2,P_3$
be three inflexion points of $C$, with $L_i$ the tangent lines at
$P_i$ to $C$.  Taking $B=C+L_1+L_2+L_3$ we construct the multiple
cyclic plane with three sheets $S_0$ associated to $B$ and
$H_\infty$.  The curve $B$ has three points of type $A_5$ at the
$P_i$'s, hence $n=3$, $b=6$ and $r=3$ in Remark \ref{r:oddEven} i).
We have 
\[
  q(S) 
  = h^1(\PP^2,\Ii_{\{P_1,P_2,P_3\}}(1)).
\]
So, if the three inflexion points are chosen on a line, then the
irregularity is $1$.  If the points are not aligned, then the
irregularity is $0$.  These two configurations give an example of a
Zariski pair.

\subsection*{Artal-Bartolo's second example in \cite{Ar}}
Let $P$ be a fixed point and $K=\{P_1,\ldots,P_9\}$ a cluster centered
at $P$, all its points being free.  It represents  a curvilinear
subscheme $Z=Z_K$.  In \cite{Ar}, Artal-Bartolo  considers sextics
with an $A_{17}$ type singularity at $P$, with $P_2,\ldots,P_9$ the
infinitely near points of the minimal resolution. 

1) If $P_3$ lies on the line $L$ determined by $P_1$ and $P_2$ and if
$K$ does not impose independent conditions on cubics, then all sextics
are reducible.  Let $B$ be the union of two smooth cubics from
$|\Ii_Z(3)|$, and let $H_\infty$ be a line transverse to $B$.  
If $S_0$ is the $3$-cyclic multiple plane associated to $B$ and
$H_\infty$, then by Remark \ref{r:oddEven} i), 
\[
  q(S) = h^1(\PP^2,\Ii_{Z^{[3]}}(1)) = 1.
\]
Similarly, if $S_0$ is the $6$-cyclic multiple plane, then 
\[
  q(S) = h^1(\PP^2,\Ii_{Z^{[3]}}(1))
  +h^1(\PP^2,\Ii_{Z^{[6]}}(2)) 
  = 2,
\]
since there is no irreducible conic through $Z^{[6]}$---\ie through
the points $P_1,\ldots,P_6$--- but the double line 
$2L$: if $K'=\{P_1^2,P_2^2,P_3^2\}$, then $Z^{[6]}\subset Z_{K'}$. 

More generally, if $S_0$ is the $n$-cyclic multiple plane associated
to $B$ and $H_\infty$, then by the same argument it follows that 
$q(S)=2$ when $n\equiv0\mod 6$, $q(S)=1$ when $n\equiv3\mod 6$, and
$q(S)=0$ otherwise. 

2) If $P_3\notin L$ and $P_6\in\Gamma$, the conic through
$P_1,\ldots,P_5$, then there exists an irreducible sextic  with an
$A_{17}$ type singularity at $P$, such that the  intersection with
$\Gamma$ is supported only at $P$.   
If $S_0$ is the $n$-cyclic multiple plane associated to $B$ and to a
transverse line to it, then 
\[
  q(S) = h^1(\PP^2,\Ii_{Z^{[6]}}(2)) = 1
\]
when $n$ is divisible by $6$, and $q(S)=0$ otherwise.

3) If $P_3\notin L$ and $P_6\notin\Gamma$, the conic through 
$P_1,\ldots,P_5$, then, for every reduced sextic $B$ with an 
$A_{17}$ type singularity at $P$, if $S_0$ is the $n$-cyclic 
multiple plane associated to $B$ and to a transverse line 
to it, then $q(S)=0$.

\begin{rem*}
In \cite{Ar} it is shown that in the third case above, two 
configurations may appear: either $P_1,\ldots,P_9$ do not impose 
independent conditions on cubics and $B$ is the union of two 
smooth cubics, or the points impose independent conditions on cubics
and $B$ is irreducible.  Using these and the two 
configurations in 1) and 2), two more Zariski couples are thus
produced there. 
\end{rem*}

\subsection*{Oka's example in \cite{Ok}}
In \cite{Ok}, when $p$ and $q$ are relatively prime integers, Oka
constructs the curve $C_{p,q}$ of degree $pq$ enjoying the following
property: $C_{p,q}$ has $pq$ cusp singularities each of which is
locally defined by the equation $x^p+y^q=0$.  We shall show that the
$pq$-multiple plane associated to $C_{p,q}$ is irregular, the
irregularity being equal to $(p-1)(q-1)/2$.  

We start with the particular case $p=2$, since all ideas of the
general computation are already present in this situation.  The
construction of the branching curve $B=C_{2,2m+1}$ is as follows.  Let
$C\subset\PP^2$ be a curve of degree $2m+1$ and let $\Gamma$ be a
conic transverse to $C$.  If $f=0$ and $g=0$ are homogeneous equations
for $C$ and $\Gamma$ respectively, then the curve $B$ is defined by
$f^2+g^{2m+1}=0$.  It is a curve of degree $4m+2$ with $4m+2$ singular
points of type $A_{2m}$.  Let $S_0$ be the $(4m+2)$-cyclic multiple
plane associated to $B$ and let $S$ be the normal cyclic covering
constructed in Section \ref{s:theIrregularity}.

\paragraph{Claim} $q(S)=m$.\\
To see this, we apply Remark \ref{r:oddEven} ii) to obtain
$q(S)=\sum_{\alpha=1}^m 
  h^1(\PP^2,\Ii_{\Zz^{[\alpha]}}(2m+2\alpha-2))$, 
where $\Zz^{[\alpha]}=\cup_P Z_P^{[\alpha]}$ and $Z_P^{[\alpha]}$ is
the curvilinear subscheme associated to the cluster
$\{P_1=P,P_2,\ldots,P_{m+2}\}$.  We shall show that all the terms of
the sum equal $1$.  To do this, we apply the trace-residual exact
sequence with respect to $\Gamma$, see \cite{Hi} or Remark
\ref{r:trRes}, and obtain
\[
  0 
  \lra \Ii_{\Zz^{[\alpha-1]}}(2m+2\alpha-4)
  \lra \Ii_{\Zz^{[\alpha]}}(2m+2\alpha-2)
  \lra \Oo_{\PP^1}(4\alpha-6)
  \lra 0.
\]
Since $C\in|\Ii_{\Zz^{[m+1]}}(2m+1)|$, the map
$H^0(\PP^2,\Ii_{\Zz^{[\alpha]}}(2m+2\alpha-2))\to 
  H^0(\PP^1,\Oo_{\PP^1}(4\alpha-6))$ 
from the long exact sequence in cohomology is surjective for every
$1\leq\alpha\leq m$.  Hence 
\[
  h^1(\PP^2,\Ii_{\Zz^{[r]}}(4r-2))
  = \cdots = h^1(\PP^2,\Ii_{\Zz^{[1]}}(2r))
  = h^1(\PP^1,\Oo_{\PP^1}(-2))
  = 1.
\]

\begin{rem*}
The irregularity of the $n$-cyclic multiple plane associated to
$B$ and to a line $H_\infty$ transversal to $B$, $n$ being an
arbitrary positive integer, may be computed by the same argument.  
Of course, if $2m+1$ is a prime number, then $q(S)=0$ unless $4m+2$
divides $n$, see Corollary \ref{c:qIrr}.  But if $2m+1$ is not a prime
number, then irregular cyclic multiple planes exist for other values
of $n$. For example, if $2m+1=15$ and $n=40$, then 
\[
  q(S)
  = h^1(\PP^2,\Ii_{\Zz^{[3]}}(18))+h^1(\PP^2,\Ii_{\Zz^{[6]}}(24))
  = 2.
\]
\end{rem*}

In the general case, if $p<q$, let $B=C_{p,q}$ and $C_p$ and $C_q$ be
the smooth curves of degree $p$ and respectively $q$ used in the
construction of $B$.  $C_p$ and $C_q$ intersect transversely and if
$P$ is an intersection point, then $B$ has a singularity at $P$ given
locally by an $x^p+y^q=0$ type equation. 

\paragraph{Claim} $q(S)=(p-1)(q-1)/2$.\\
By Corollary \ref{c:pqTypeSing},
\[
  q(S) 
  = \sum_{\substack{\alpha,\beta\geq 1 \\ \alpha p+\beta q<pq}}
  h^1(\PP^2,\Ii_{\Zz_{p,q}(\widetilde{\alpha,\beta})}
  (-3+\alpha p+\beta q)).
\]
The sum consists of $(p-1)(q-1)/2$ terms, and as before, we shall show
that each of them equals $1$.  For an arbitrary couple
$(\alpha,\beta)$, with $\alpha\geq2$, we first apply the
trace-residual exact sequence $\alpha-1$ times with respect to $C_p$.
We have
\[
  0 
  \to \Ii_{\Zz(\widetilde{\alpha-1,\beta})}(-3+(\alpha-1)p+\beta q)
  \to \Ii_{\Zz(\widetilde{\alpha,\beta})}(-3+\alpha p+\beta q)
  \stackrel{\rho}{\to}
  \Ii_{\Trace_{C_p}\Zz(\widetilde{\alpha,\beta})}(-3+\alpha p+\beta q)
  \to 0.
\]
Let $w_1$ be the weight of $P_1=P$ in the cluster
$K_{p,q}(\widetilde{\alpha,\beta})$.  Since the cluster is not greater
than $K_{p,q}(0,\floor{(\alpha-1)p/q}+\beta)$, it is easy to see that 
$w_1\leq\floor{(\alpha-1)p/q}+\beta$.
Then $\Zz(\widetilde{\alpha,\beta})\subset w_1C_q$ and together with
the identity
\[
  -3+\alpha p+\beta q 
  = -3+p+\left((\alpha-1)p-\floor{\frac{(\alpha-1)p}{q}}q\right)+
  \left(\floor{\frac{(\alpha-1)p}{q}}+\beta-w_1\right)q+w_1q
\]
imply the surjectivity of $H^0\rho$.  We conclude that 
\begin{equation}
  \label{eq:onCp}
  h^1(\PP^2,\Ii_{\Zz(\widetilde{\alpha,\beta})}(-3+\alpha p+\beta q))=
  h^1(\PP^2,\Ii_{\Zz(\widetilde{1,\beta})}(-3+p+\beta q)) 
\end{equation}
whenever $\alpha\geq 2$.  Then, in case $\beta\geq2$, we apply
$\beta-1$ times the trace-residual exact sequence with respect to
$C_q$ starting with the subscheme $\Zz(\widetilde{1,\beta})$.  As
before, we have
\[
  0 
  \to \Ii_{\Zz(\widetilde{1,\beta-1})}(-3+p+(\beta-1)q)
  \to \Ii_{\Zz(\widetilde{1,\beta})}(-3+p+\beta q)
  \stackrel{\rho}{\to}
  \Ii_{\Trace_{C_q}\Zz(\widetilde{1,\beta})}(-3+p+\beta q)
  \to 0,
\]
the surjectivity of $H^0\rho$ being given by the inequality
$w<1+(\beta-1)q/p$, with $w$ the sum of the weights of the points
$P_{1,1},\ldots,P_{1,a_1},P_{2,1}$ in $K_{p,q}(\widetilde{1,\beta})$
and the inclusion $\Zz(\widetilde{1,\beta})\subset wC_p$.  So
\begin{equation}
  \label{eq:onCq}
  h^1(\PP^2,\Ii_{\Zz(\widetilde{1,\beta})}(-3+p+\beta q)) =
  h^1(\PP^2,\Ii_{\Zz(\widetilde{1,1})}(-3+p+q)).
\end{equation}
Finally, since $\Zz_{p,q}(\widetilde{1,1})=\cup_PP$, we apply once
more the trace-residual exact sequence of $\Zz(\widetilde{1,1})$ 
with respect to $C_p$ and get
\[ 
  0 
  \lra \Oo_{\PP^2}(-3+q)
  \lra \Ii_{\Zz(\widetilde{1,1})}(-3+p+q)
  \lra \Oo_{C_p}(-3+p)
  \lra 0.
\]
Since $q>p$, 
$h^1(\PP^2,\Ii_{\Zz(\widetilde{1,1})}(-3+p+q))=h^1(C_p,\Oo_{C_p}(-3+p))=1$.
Together with (\ref{eq:onCp}) and (\ref{eq:onCq}) finish the claim.

\begin{rem}[The trace-residual exact sequence]
  \label{r:trRes}
Let $X$ be a projective variety, $D$ be a Cartier divisor on $X$ and
$Z$ be a a closed subscheme of $X$. 
The schematic intersection $\Trace_DZ=D\cap Z$ defined by the ideal
sheaf $(\Ii_D+\Ii_Z)/\Ii_D$ is called the trace of $Z$ on $D$.  The
closed subscheme $\Res_DZ\subset X$ defined by the conductor ideal
$(\Ii_Z:\Ii_D)$ is called the residual of $Z$ with respect to $D$.
Following \cite{Hi}, the canonical exact sequence 
\begin{equation*}
  0\lra
  \Ii_{\Res Z}(-D)\lra 
  \Ii_Z\lra 
  \Ii_{\Trace Z}\lra 0.
\end{equation*}
is called the trace-residual exact sequence of $Z$ with respect
to $D$. 
\end{rem}

\subsection*{A specialization of Oka's example when $p=2$}

Keeping the notation from the first part of the previous paragraph,
the conic $\Gamma$ is now the union of two distinct lines that
intersect at $O$ and $C$ is a smooth curve of degree $2m+1$ passing
through $O$ and intersecting transversely the lines of $\Gamma$ at
this point.  The curve $B$ has $4m$ points of type $A_{2m}$ and one
singular point at $O$ of type $A_{4m+1}$.  It can be shown that the
irregularity of the $(4m+2)$-cyclic multiple plane associated to $B$
is again $m$.  We develop the computation for $m=2$.  In this case,
$B$ is a curve of degree $10$ with $8$ points of type $A_4$ and one
point of type $A_9$.  By Theorem \ref{th:q} and using the notation
from Remark \ref{r:oddEven}, the irregularity is given by
\[
  h^1(\PP^2,\Ii_{\xi^{[1]}\cup Z_O^{[2]}}(4)) 
  + h^1(\PP^2,\Ii_{{\xi^{[2]}}\cup Z_O^{[4]}}(6)),
\]
where $\xi^{[1]}$ is the support of the points of type $A_4$ and
$\xi^{[2]}=\cup_{P\text{ of type }A_4}Z_P^{[2]}$ is the support plus
the tangent directions.  Now, $10$ points on a conic do not impose
independent conditions on quartics, hence the first term is $1$.  The
second term is seen to be equal to the first after applying the
trace-residual exact sequence with respect to the two lines of
$\Gamma$.  So the irregularity is $2$.

The computations for $m=1$ lead to a branching curve of
degree $6$ with $4$ cusps and an $A_5$ singularity at $O$.  The
irregularity of a $6$-cyclic multiple plane is $1$, given
by $h^1(\PP^2,\Ii_{\xi^{[1]}\cup Z_O^{[2]}}(2))$. 
If in addition, the two lines of the degenerate conic $\Gamma$ are
brought together such that the cusps collapse two by two, the
branching curve has $3$ $A_5$ singularities.  For a $6$-multiple
plane, $q=2$, with the contributions of the superabundance of the
singularities with respect to the lines and the conics both equal to
$1$.  Necessarily, by Corollary \ref{c:qIrr}, the branching curve is
reducible; it is Artal-Bartolo's first example.

\subsection*{Line arrangements following \cite{Es}}

In this example we consider as branch curve a line arrangement
$B=\cup_{i=1}^bL_i\subset\PP^2$ that has only nodes and ordinary
triple points as singularities.  For an ordinary triple point $2/3$ is
the only subunitary jumping number.  
By Corollary \ref{c:pqTypeSing}, 
if $H_\infty$ is a line transverse to $B=\cup_{i=1}^bL_i$, then the
normal $n$-cyclic covering $S$ corresponding to the $n$-cyclic
multiple plane associated to $B$ and $H_\infty$ is irregular if and
only if $3$ divides both $b$ and $n$, and $|\Ii_\Zz(-3+\frac{2b}{3})|$
is superabundant, in which case 
\[
  q(S) = h^1(\PP^2, \Ii_\Zz(-3+\frac{2b}{3})).
\]
In case $S$ is irregular, it can be shown that the irregularity is
bounded by a constant depending on the arrangement $B$.  More
precisely, we have 

\begin{pro} \label{p:lineArrBound}
Let $B=\cup_{i=1}^b L_i$, $H_\infty$ and $S$ be as above with $b$ and
$n$ divisible by $3$.  
If $t_i$ is the number of triple points lying on the line $L_i$ for
each $i$, then 
\[
  q(S) \leq \min_{i=1}^b t_i.
\]
\end{pro}

For the proof (see \cite{Es} for a different argument), we will start
with a preliminary lemma.

\begin{lem}
If $3$ divides both $b$ and $n$ and if one line of the arrangement
contains no triple point, then $q(S)=0$. 
\end{lem}
\proof
Let $B'$ be the arrangement of the $b-1$ lines of $B$ except the one
one with no triple point.  
If $S'$ is the normal $n$-cyclic covering corresponding to the
$n$-cyclic multiple plane associated to $B'$ and $H_\infty$, then
$q(S')=0$ since $3$ does not divide $\deg B'$.  Taking $k=2n/3$ and
denoting by $\Zz$ the support of the triple points, we obtain 
\[
  0 = h^1(\PP^2,\Ii_\Zz(-3+\ceil{\frac{2(b-1)}{3}})
  = h^1(\PP^2,\Ii_\Zz(-3+\frac{2b}{3}) = q(S).
\] 
\qed

\likeproof[Proof of Proposition \ref{p:lineArrBound}]
Let us suppose that $L_1$ is the line containing the minimum number of
triple points.  If $B'=L'_1\cup\bigcup_{i\neq1}L_i$ is a line
arrangement with  no triple point on $L'_1$, then by the previous lemma,
$h^1(\PP^2,\Ii_{\Res_{L_1}\Zz}(-3+2b/3))=0$.  But
\[
  h^1(\PP^2,\Ii_\Zz(-3+\frac{2b}{3})) 
  \leq h^1(\PP^2,\Ii_{\Res_{L_1}\Zz}(-3+\frac{2b}{3})) 
  + \card(\Zz-\Res_{L_1}\Zz) = t_1,
\] 
hence the result.
\qed

\begin{exa*}
Let $B$ be the line arrangement of $9$ lines with $9$ triple points
represented below.  In a convenient affine coordinate system $(x,y)$,
the triple points that lie in the affine plane are the following: 
$(0,0),(\pm2,-2),(-2,0),(0,s),(2,s)$ and $2s/(s+4)\,(-1,1)$, with 
$s\neq-2,0$ and $2$. 
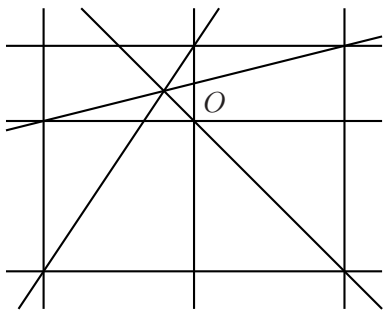
\begin{center}
  \begin{pspicture}(-2.25,-2.25)(1.75,1.75)%\showgrid
    \psset{unit=1cm}
    \psline(-2,-2.5)(-2,1.5)
    \psline(0,-2.5)(0,1.5)
    \psline(2,-2.5)(2,1.5)
    \psline(-2.5,1)(2.5,1)    
    \psline(-2.5,0)(2.5,0)
    \psline(-2.5,-2)(2.5,-2)

    \psline(-2.3333,-2.5)(0.3333,1.5)
    \psline(-2.5,-.125)(2.5,1.125)
    \psline(2.5,-2.5)(-1.5,1.5)

    \uput{\labelsep}[ur](0,0){$O$}
\end{pspicture}
\end{center}
It is easy to see that there are two cubics---each the union of three
lines---through the $9$
triple points, \ie the system of cubics through the points is
superabundant.  It follows that the irregularity of the $n$-cyclic
multiple plane associated to $B$ and to a line $H_\infty$ transverse to
$B$, is $1$ if and only if $3$ divides $n$.  

If $s=2$, then the arrangement specialize to an arrangement with $10$
triple points, $4$ of them lying on the line $x+y=0$.  But these
points lie on a cubic, the union of three of the lines of $B$, and
again $h^1(\PP^2,\Ii_\Zz(3))=1$, hence the irregularity is $1$ in this
case too. 
\end{exa*}

\begin{rem}
  \label{r:nonTransverse}
The irregularity depends on the position of the line $H_\infty$ with
respect to $B$.  To see this, let $B$ be the line arragement below of
$5$ lines with $2$ triple points.
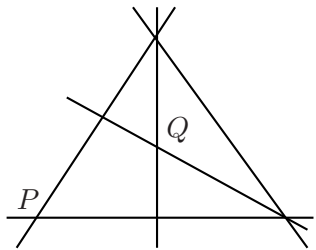
\begin{center}
  \begin{pspicture}(-2.25,-2.25)(1.75,1.75)%\showgrid
    \psset{unit=.8cm}
    \psline(0,-2.5)(0,1.5)
    \psline(-2.5,-2)(2.6,-2)
    \psline(-2.3333,-2.5)(0.3,1.5)
    \psline(2.5,-2.5)(-.4,1.5)
    \psline(2.5,-2.2)(-1.5,0)

    \uput{\labelsep}[ur](0,-.9){$Q$}
    \uput{\labelsep}[ul](-1.8,-2){$P$}
\end{pspicture}
\end{center}
If $H_\infty$ is transverse to $B$, then the irregularity of the
$6$-cyclic multiple plane is $0$.  But if $H_\infty$ is the line
through the double points $P$ and $Q$ then the irregularity jumps to
$1$.
\end{rem}

\begin{flushleft}
  Daniel Naie\\
  D\'epartement de Math\'ematiques\\
  Universit\'e d'Angers\\
  F-40045 Angers\\
  France
\end{flushleft}
\end{document}